\documentclass{m2an}
\usepackage[colorlinks=true,citecolor=blue,urlcolor=blue]{hyperref}
\usepackage{amssymb} 
\usepackage{graphicx}
\newcommand {\bx} {{\mathbf x}}
\newcommand {\rw} {{\mathrm w}}
\newcommand {\rv} {{\mathrm v}}
\DeclareMathOperator\erf{erf}
\newtheorem{remark}{Remark}[section]

\newcommand{\rev}[1]{{\color{black} #1}}
\begin{document}
\title{\rev{Projected gradient} stabilization of sharp  and diffuse interface formulations in unfitted Nitsche finite element methods}
\author{Maxim Olshanskii}
\address{Department of Mathematics, University of Houston, Houston, Texas 77204-3008, USA;
\email{maolshanskiy@uh.edu}}
\author{Jan-Phillip Bäcker}
\address{Institute of Applied Mathematics (LS III), TU Dortmund University, Vogelpothsweg 87, D-44227 Dortmund, Germany;
\email{jan-phillip.baecker@tu-dortmund.de,\ kuzmin@math.uni-dortmund.de}}
\author{Dmitri Kuzmin}
\sameaddress{2}
\date{\today}
\begin{abstract}
We introduce an unfitted Nitsche finite element method with a new ghost-penalty stabilization based on local projection of the solution gradient.
The proposed ghost-penalty operator is straightforward to implement, ensures algebraic stability, provides an implicit extension of the solution beyond the physical domain, and stabilizes the numerical method for problems dominated by transport phenomena. This paper presents both a sharp interface version of the method and an alternative diffuse interface formulation designed to avoid integration over implicitly defined embedded surfaces. A complete numerical analysis of the sharp interface version is provided.
The results of several numerical experiments support the theoretical analysis and illustrate the performance of both variants of the method.
\end{abstract}
%
%
\maketitle
\section*{Introduction}
\label{sec:intro}

During the last two decades, immersed interface finite element methods have become an increasingly popular approach to solving complex problems with discontinuous material properties and moving boundaries \cite{bordas}. The ability to perform simulations on a fixed (and often structured) mesh offers significant savings in terms of computational cost. However, the stability and accuracy of such geometrically unfitted formulations are strongly influenced by the numerical treatment of interface conditions. The~unfitted Nitsche \cite{hansbo2,wadbro} finite element method (FEM) and its stabilized counterparts, commonly known as cutFEM \cite{cutfem}, are among the most successful discretization techniques of this type. The underlying weak forms include surface integrals that penalize violations of the interface conditions, as well as ghost penalty terms \cite{burman2010,burman2012,olsh} that prevent ill-conditioning due to the presence of small cut cells. In many cases, optimal a priori error estimates and upper bounds for the condition number of the system matrix can be obtained for such unfitted finite element discretizations.

The practical implementation of cutFEM in existing codes presents certain challenges, including numerical integration over embedded and implicitly defined domains and interfaces, as well as the proper extension of the solution outside the physical domain. To address these challenges in the context of a stabilized unfitted Nitsche method, we:
\begin{itemize}
\item[a)] propose a versatile stabilization-extension operator for cutFEM approaches;
\item[b)] analyze the stability and accuracy of the resulting finite element scheme;
\item[c)] introduce and numerically study a diffuse interface version of our method.
\end{itemize}

In contrast to some other ghost penalties, our new projected-gradient (PG) stabilization is parameter-free and easy to implement \rev{for finite elements of any order}. It smoothly extends subdomain solutions beyond the physical domain into a buffer zone of a desired width and offers the additional benefit of stabilizing the convective terms, if present in the governing equation. \rev{Convective terms can also be stabilized using ghost penalties of the continuous interior penalty (CIP) type \cite{CIP}. Although CIP stabilization is effective in many contexts, it does not readily accommodate limiters for enforcing discrete maximum principles (as in \cite{WSbook,entropyHO,lohmann2017}) and offers no direct approach for constructing compact-stencil preconditioners for efficiently solving the resulting linear systems. In comparison, our PG stabilization is naturally suited to both tasks. Furthermore, we numerically demonstrate its optimal accuracy for high-order finite elements in this paper. }

The results of our theoretical studies confirm the favorable properties of our sharp interface PG method for an elliptic interface problem. In particular, we prove stability and derive a priori error estimates. \rev{The close relationship with the ghost penalty stabilization term \cite[Eq. 23]{badia} of an aggregated unfitted method is noticed and discussed.} In the diffuse interface version, surface integrals are approximated by volume integrals using a regularized delta function, as in \cite{hysing,kublik2016,kublik2018,unfitted}. The extension of data from the interface into its neighborhood is performed using Taylor expansions and the fast closest-point search algorithm developed in \cite{unfitted}.

The remainder of this paper is organized as follows. In Section~\ref{sec:unfitted}, we introduce a model problem and review the classical Nitsche finite element method. The new local PG stabilization operator is introduced in Section~\ref{sec:sharp} \rev{for elements of arbitrary order} in the sharp interface context. The diffuse interface formulation, presented in Section~\ref{sec:diffuse}, eliminates the need for integration over sharp interfaces and cut cells. In Section~\ref{sec:analysis}, we analyze the PG stabilization technique \rev{for the particular case of affine elements}. Numerical results are presented in Section~\ref{sec:num}. The convergence rates for test problems with smooth exact solutions are shown to be optimal. When the gradient is discontinuous across the interface, kinks in the solution profile are well captured. Finally, the main findings are summarized in Section~\ref{sec:conclusions}.

\section{An interface problem and unfitted Nitsche FE method}
\label{sec:unfitted}

Let $\Omega\subset\xR^d$, $d\in\{2,3\}$ be a bounded Lipschitz domain, which is divided into two subdomains, $\Omega_1$ and $\Omega_2$, by a smooth embedded interface~$\Gamma$. 
For a function $v:\Omega\to\mathbb{R}$, its  restriction  to $\Omega_k,\ k=1,2$
is denoted by $v_k$. 
The jump of $v$ across $\Gamma$ is denoted by $[\![ v]\!]=v_1-v_2$. An average $\{v\}=\kappa_1v_1+\kappa_2v_2$  is defined using a pair of weights 
$\kappa_1\ge 0$ and $\kappa_2\ge 0$ such that $\kappa_1+\kappa_2=1$. The values of $\kappa_1$ and $\kappa_2$
may depend on $\mathbf{x}\in\Gamma$. The jump and average of a vector field $\mathbf{v}(\mathbf{x})$ are defined by
$[\![ \mathbf v]\!]=(\mathbf{v}_1-\mathbf{v}_2)\cdot\mathbf{n}$ and $\{\mathbf v\}=(\kappa_1\mathbf v_1+\kappa_2\mathbf v_2)\cdot
\mathbf{n}$, respectively.
\smallskip

We are interested in  the following elliptic interface problem
\begin{subequations}\label{eq:inter}
\begin{align}
	-\nabla\cdot(\mu\nabla u)&=f\quad\text{in }\Omega_k,~{k=1,2},\label{eq:pde}\\
	u &= 0\quad\text{on }\partial\Omega,\label{eq:bc1}\\
	[\![ u]\!]&=0\quad\text{on }\Gamma,\label{eq:bc2}\\
	[\![ \mu\nabla u]\!]&=0\quad\text{on }\Gamma
	\label{eq:bc3},
\end{align}
\end{subequations}
where
$$
\mu=\begin{cases}
\mu_1 & \mbox{in}\ \Omega_1,\\
\mu_2 & \mbox{in}\ \Omega_2,
\end{cases}
$$
$\mu_k>0$ are constant rates of diffusive transport, and $f\in L^2(\Omega)$.

To build a numerical method for solving \eqref{eq:inter}, we consider a regular conforming triangulation  $\mathcal{T}_h$ of $\Omega$. For simplicity, we assume that the triangulation is fitted to the domain $\Omega$ in the sense that  $\bigcup_{K\in\mathcal{T}_h}K
=\bar\Omega$. For any $\delta\ge0$, we define the sub-triangulations
\[
\mathcal{T}_{h,k}(\delta)=\{K\in \mathcal{T}_h\,:\, \exists\, \mathbf{x}\in K\ \mbox{s.t.}\ \mbox{dist}(\bx,\Omega_k)\le \delta\},\quad k=1,2,
\]
consisting of simplices that intersect  a
$\delta$-neighborhood of $\Omega_k$. The notation
$$
\Omega_{h,k}(\delta)=\mathrm{int}\big(\bigcup_{K\in \mathcal{T}_{h,k}(\delta)} K\big),\quad k=1,2
$$
is used for the corresponding mesh-dependent domains. Note that $\Omega_k\subset \Omega_{h,k}(\delta)$ for all $\delta\ge0$. In particular, for $\delta=0$ the closure of $\Omega_{h,k}(0)$ consists of simplices that have a non-empty intersection with $\Omega_k$.

Denote by $V_{h,k}(\delta)$ the finite element spaces of continuous piecewise \rev{polynomial functions of degree $p>0$} with respect to $\mathcal{T}_{h,k}(\delta)$:
\[
V_{h,k}(\delta) =\{v_h\in C(\Omega_{h,k}(\delta)):\, v_h|_K\in P_{p}(K)~\forall\,K\in \mathcal{T}_{h,k}(\delta)~\text{and}~ v_h=0~\text{on}\,\partial\Omega\cap\partial\Omega_{h,k} \}.
\]

The standard  Nitsche finite element method for  \eqref{eq:inter} reads as follows~\cite{hansbo2,nitsche}: Find $u_{h,k}\in V_{h,k}(0)$, $k=1,2$, such that 
\begin{align}
\begin{split}
	\int_{\Omega_1}\mu_1\nabla u_{h,1}\cdot\nabla w_{h,1}\xdif \bx&+\int_{\Omega_2}\mu_2\nabla u_{h,2}\cdot\nabla w_{h,2}\xdif \bx\\&-\int_\Gamma[\![u_h]\!]\{\mu\nabla w_h\}\xdif s-\int_\Gamma\{\mu\nabla u_h\}[\![w_h]\!]\xdif s+\int_\Gamma\alpha[\![u_h]\!][\![w_h]\!]\xdif s\\=&\int_{\Omega_1}f_1w_{h,1}\xdif \bx+\int_{\Omega_2}f_2w_{h,2}\xdif \bx\quad \forall\,w_{h,k}\in V_{h,k}(0),~k=1,2.
\end{split}\label{nitschewf}
\end{align}
Following~ \cite{hansbo2}, we define the averages using weights $\kappa_1,\kappa_2$ that are constant in each cell and represent the
volume fractions of cut cells $K$ from
$\mathcal T_h$: 
\begin{align}
\kappa_k|_{K}=\frac{|K\cap\Omega_k|}{|K|},\qquad k=1,2.
\end{align}
The interior penalty parameter $\alpha$ depends on the local mesh size $h_K=\mathrm{diam}(K)$: $\alpha|_{K}= \alpha_0h^{-1}_K$ with a sufficiently large $\alpha_0=O(1)$.

We define the bilinear and linear forms: 
\begin{align*}
a_h(u,w)&:=\int_{\Omega_1}\mu_1\nabla u_{1}\cdot\nabla w_{1}\xdif \bx+\int_{\Omega_2}\mu_2\nabla u_{2}\cdot\nabla w_{2}\xdif \bx\\&\qquad -\int_\Gamma[\![u]\!]\{\mu\nabla w\}\xdif s-\int_\Gamma\{\mu\nabla u\}[\![w]\!]\xdif s+\int_\Gamma\alpha[\![u]\!][\![w]\!]\xdif s,\\
f(w)&:=\int_{\Omega_1}f_1w_{1}\xdif \bx+\int_{\Omega_2}f_2w_{2}\xdif \bx,\qquad u,w\in H^1(\Omega_1)\times H^1(\Omega_2)
\end{align*}
to write \eqref{nitschewf} as the abstract problem
\begin{align}
a(u_h,w_h) = f(w_h)\quad \forall\,w_h\in V_{h,1}(0)\times V_{h,2}(0).
\label{nitschewf2}
\end{align}

Recalling the definitions of jumps and averages, we find the decomposition
$$
[\![u_h]\!]\{\mu\nabla w_h\}+\{\mu\nabla u_h\}[\![w_h]\!]
-\alpha[\![u_h]\!][\![w_h]\!]=q_1(u_{h,1},u_{h,2},w_{h,1})-
q_2(u_{h,1},u_{h,2},w_{h,2}),
$$
where
\begin{subequations}\label{g12inter}
\begin{align}
	q_1(u_1,u_2,w_1)&=(u_1-u_2)(\kappa_1\mu_1\nabla w_1)\cdot\mathbf{n}
	+(\kappa_1\mu_1\nabla u_1+\kappa_2\mu_2\nabla u_2)w_1\cdot\mathbf{n}
	-\alpha(u_1-u_2)w_1,\\
	q_2(u_1,u_2,w_2)&=(u_2-u_1)(\kappa_2\mu_2\nabla w_2)\cdot\mathbf{n}
	+(\kappa_1\mu_1\nabla u_1+\kappa_2\mu_2\nabla u_2)w_2\cdot\mathbf{n}
	+\alpha(u_2-u_1)w_2.
\end{align}
\end{subequations} 
Using this notation, the finite element method \eqref{nitschewf} can be also written as the system of two coupled problems:
Find $u_{h,k}\in V_{h,k}(0)$, $k=1,2$ such that 
\begin{subequations}\label{subq1q2}
\begin{align}
	\int_{\Omega_1}\mu_1\nabla u_{h,1}\cdot\nabla w_{h,1}\xdif \bx&
	-\int_\Gamma q_1(u_{h,1},u_{h,2},w_{h,1})\xdif s
	=\int_{\Omega_1}f_1w_{h,1}\xdif \bx\quad \forall\,w_{h,1}\in V_{h,1}(0),\\
	\int_{\Omega_2}\mu_2\nabla u_{h,2}\cdot\nabla w_{h,2}\xdif \bx&+\int_\Gamma
	q_2(u_{h,1},u_{h,2},w_{h,2})
	\xdif s=\int_{\Omega_2}f_2w_{h,2}\xdif \bx\quad \forall\,w_{h,2}\in V_{h,2}(0).
\end{align}
\end{subequations}
In what follows, we consider stabilized versions of \eqref{nitschewf}. They can be restricted to the subdomains $\Omega_{h,1}(\delta)$ and $\Omega_{h,2}(\delta)$ similarly by using test functions $w_{h,1}\in V_{h,1}(\delta)$ and $w_{h,2}\in V_{h,2}(\delta)$, respectively.

\section{Sharp interface PG formulation}
\label{sec:sharp}

Although well-posed, the finite element (FE) formulation \eqref{nitschewf} is prone to numerical instabilities. These instabilities arise from the potential inclusion of simplices in $\mathcal{T}_{h,k}(0)$ that intersect the `physical' domain $\Omega_k$ only minimally. This stability issue was first addressed in \cite{burman2010} by introducing ghost penalty stabilization on simplices intersected by and adjacent to $\Gamma$. The approach proposed in \cite{olsh} extends this penalization to a broader strip of elements within a distance $\delta > 0$ from $\Gamma$, as part of an implicit extension procedure for an Eulerian FE discretization. This formulation is particularly relevant to problems where the evolving interface $\Gamma$ represents the boundary of a time-dependent embedded domain.

Let us now introduce a gradient-projection stabilization technique as an alternative to traditional ghost penalties in both roles: stabilizing against small cuts and providing a proper extension of $u_{h,k}$ to $\Omega_{h,k}(\delta)$ for $\delta > 0$. We define 
the \emph{\rev{\rev{projected gradient}}} (PG) stabilization term
\begin{align}\label{neumannpen}
s_{h}(u,w):=
\int_{\Omega_{h,1}(\delta)}\mu_1(\nabla u_{1}-\mathbf{g}_{h,1})\cdot\nabla w_{1}\xdif \bx
+\int_{\Omega_{h,2}(\delta)}\mu_2(\nabla u_{2}-\mathbf{g}_{h,2})\cdot\nabla w_{2}\xdif \bx
\end{align}
with $ u,w\in H^1(\Omega_{h,1}(\delta))\times H^1(\Omega_{h,2}(\delta))$.
The stabilization uses  lumped-mass $L^2$
projections $\mathbf{g}_{h,k}$ of the (generally discontinuous) gradients $\nabla u_{k}
\in (L^2(\Omega_{h,k}(\delta)))^d$ into the spaces $(V_{h,k}(\delta))^d$. Those $\mathbf{g}_{h,k}$ are computed locally as follows:
Consider  the nodal basis  $\{\varphi_j\}_{j=1,\dots,N_{h,k}}$ in $(V_{h,k}(\delta))$,  $N_{h,k}=\mathrm{dim}(V_{h,k}(\delta))$ and let
\begin{equation}\label{aux204}
\mathbf{g}_{h,k}=\sum_{j=1}^{N_{h,k}}\mathbf{g}_{j,k}\varphi_j,\qquad k=1,2.
\end{equation}

\begin{itemize}
\item{\rev{\bf Affine elements.}} 
For the lowest order case of linear element ($p=1$), the coefficients  $\mathbf{g}_{j,k}$ are explicitly given by
\begin{align}\label{lumpedg}
\mathbf{g}_{j,k}=\frac{\int_{\Omega_{h,k}(\delta)}\varphi_j \nabla u_{k}\xdif \bx}{\int_{\Omega_{h,k}(\delta)}\varphi_j\xdif \bx},\qquad k=1,2.
\end{align}
Note that computing $\mathbf{g}_{j,k}$ involves integration only over triangles sharing the node $\bx_j$, where $\varphi(\bx_j)=1$.
\rev{The stabilization has a particularly simple algebraic form, which requires only additional gradient and lumped-mass matrices, cf. Remarks~\ref{rem22} and~\ref{rem41}.} 
\\[-1ex] 

{\color{black}
\item{\bf Higher order elements.} 
  Following~\cite{entropyHO} we define the coefficients $\mathbf{g}_{j,k}$ for
  arbitrary order finite elements as weighted averages of the one-sided limits
$\nabla u_{h,k}|_{K}(\mathbf{x}_j)$  in elements $K$ that share the point $\mathbf x_j$. Introducing the notation $\mathcal{E}_j$ for the set of such elements, we set
\begin{equation}\label{aux221g}
  \mathbf{g}_{j,k} = \frac{1}{\sum_{K \in \mathcal{E}_j}|K|}
  \sum_{K \in \mathcal{E}_j}|K|\,\nabla u_{h,k}|_{K}(\mathbf{x}_j).
\end{equation}
The definitions \eqref{lumpedg} and \eqref{aux221g} are equivalent
for affine elements.
}
\end{itemize}

Adding the PG stabilization term \eqref{neumannpen} to the finite element formulation \eqref{nitschewf2} results in the following
sharp interface formulation:  Find $u_{h,k}\in V_{h,k}(\delta)$, $k=1,2$, such that
\begin{align}
a(u_h,w_h)	+ s_{h}(u_h,w_h)  = f(w_h)\quad \forall\,w_h\in V_{h,1}(\delta)\times V_{h,2}(\delta).\label{wsiab}
\end{align}

\begin{remark}
Projection-based stabilization terms of the form \eqref{neumannpen} are widely used in FE formulations for convection-dominated transport problems \cite{codina1,codina2,john1,john2}. Consequently, their stabilizing effect can be naturally exploited in extensions of the proposed method to convection-diffusion equations (see Section \ref{sec:convection} below for a numerical example). This is one of the reasons why we apply stabilization globally on $\Omega_{h,2}(\delta)$ rather than locally in a narrow strip around the interface. Another reason is the ease of implementation. While the ghost penalty term in unfitted FEMs is typically defined within a narrow band, we are not the first to consider a global definition of the ghost penalty; see, for example, \cite{gp1,gp2}. In particular, the recent report \cite{gp3} demonstrated that defining the ghost penalty globally is necessary to ensure specific stability of the method in the case of moving interfaces.
\end{remark}

\rev{
	\begin{remark}[Algebraic representation]\label{rem22}\rm 
		Let us take a look at the algebraic form of the stabilization. To this end, using the standard Lagrangian nodal basis in our finite element spaces, we consider the matrix $L$ corresponding to the bilinear form
		$
		(\nabla \cdot, \nabla \cdot)_{L^2(\Omega_{h,1})} + (\nabla \cdot, \nabla \cdot)_{L^2(\Omega_{h,2})}.
		$
		That is, $L$ is a block-diagonal matrix consisting of two standard finite element Laplacians with mixed boundary conditions on the domains $\Omega_{h,k}$.
		Similarly, $B$ denotes the gradient matrix associated with the bilinear form
		\[
		(\nabla v_h, \mathbf{w}_h)_{L^2(\Omega_{h,1})} + (\nabla v_h, \mathbf{w}_h)_{L^2(\Omega_{h,2})},
		\]
		with $v_h \in V_{h,1} \times V_{h,2}$, $\mathbf{w}_h \in V_{h,1}^2 \times V_{h,2}^2$.
		Furthermore, $\widetilde{M}$ denotes the lumped mass matrix corresponding to the finite element space $V_{h,1}^2 \times V_{h,2}^2$.
		
		Then the matrix of the  stabilization form $s_h(\cdot,\cdot)$ in the case of linear elements ($p=1$) can be found to be
		\[
		L - B^T \widetilde{M}^{-1} B.
		\]
		
		We observe that the stabilization at the algebraic level represents the difference between $L$ and the ``mixed'' finite element discretization $B^T \widetilde{M}^{-1} B$ of the Laplacian operator.
	
	\end{remark}
}

\section{Diffuse interface PG formulation}
\label{sec:diffuse}

The calculation of  integrals over cat elements and embedded interfaces in unfitted finite element discretizations
of multidimensional interface problems is not always straightforward. To avoid the associated difficulties,  we use the methodology developed in \cite{unfitted}. 
Assume that  position of $\Gamma\subset\Omega$ is determined by a continuous level set function $\phi\in C(\bar\Omega)$, which is positive in $\Omega_1$ and negative in $\Omega_2$ 
so that
$\Gamma=\{\mathbf{x}\in\Omega:\phi(\mathbf{x}) = 0\}.
$
We assume that $\phi$ is smooth in an $O(1)$ neighborhood $\mathcal{O}(\Gamma)$ and $|\nabla\phi|\ge c>0$ in $\mathcal{O}(\Gamma)$.
Then the  vector field $\mathbf{n}=-\frac{\nabla\phi}{|\nabla\phi|}$ defines an extended unit normal to $\Gamma$, pointing outward from $\Omega_1$.
The methods employed in \cite{hysing,kublik2016,kublik2018,unfitted,zahedi} approximate surface integrals over the interface $\Gamma$
by volume integrals depending on a smoothed Dirac delta function $\delta_\epsilon(\phi)$. The small parameter $\epsilon>0$ determines the quality of the approximation to $\delta(\phi)=\lim_{\epsilon\searrow 0}\delta_\epsilon(\phi)$. A typical definition of $\delta_\epsilon(\phi)\approx\delta(\phi)$ is given below.	

For any $q\in C(\Omega)$, the integral $\int_\Gamma q(\mathbf{x})\xdif s$ can be approximated by $\int_{\Omega}q(\mathbf{x})\delta_\epsilon(\phi)|\nabla\phi|\xdif \bx$. 
As explained in \cite{unfitted}, the calculation of volume integrals that contain $\delta_\epsilon(\phi)$ requires an extension of the quantities defined on $\Gamma$ into the interior of the two subdomains. To maintain the consistency of the weak form, we extend the interface values of functions and coefficients constantly in the normal direction determined by the gradient of the level set function $\phi$. The efficient algorithm that we developed for this purpose can be found in \cite{unfitted}. Let $\mathcal{E}_{\mathrm{cp}}$ denote the closest point extension  operator such that $\mathcal{E}_{\mathrm{cp}} v(\mathbf{x})
=v(\mathbf{x}_\Gamma(\mathbf{x}))$, where $\mathbf{x}_\Gamma$ is the closest point on
$\Gamma$.
Such extension of \eqref{g12inter}
is
consistent with the jump conditions \eqref{eq:bc2}
and \eqref{eq:bc3} of the 
interface problem in the sense that
$$
\mathcal{E}_{\mathrm{cp}} q_1(u_1,u_2,w_1)-\mathcal{E}_{\mathrm{cp}} q_2(u_1,u_2,w_2)=0\qquad
\mbox{at}\ \mathbf{x}\in\Omega$$
if
$$[\![ u]\!]=[\![
\mu\nabla u]\!]=
[\![ w]\!]=0\qquad \mbox{at}\ \mathbf{x}_\Gamma(\mathbf{x})\in\Gamma.
$$
Consequently, the contribution of Nitsche integrals to the weak
form vanishes if the jumps do.
\smallskip

A diffuse interface counterpart of the sharp interface formulation
\eqref{wsiab}
is given by (cf. \cite{unfitted})
\begin{equation}\label{wdiab}
	\begin{split}
		\int_{\Omega_{h,1}(\delta)}\mu_1((1+H(\phi))\nabla u_{h,1}&-\mathbf{g}_{h,1})\cdot\nabla w_{h,1}\xdif \bx  -
		\int_{\Omega}\mathcal{E}_{\mathrm{cp}} q_1(u_{h,1},u_{h,2},w_{h,1})\delta_\epsilon(\phi)|\nabla\phi|\xdif \bx\notag \\
		+\int_{\Omega_{h,2}(\delta)}\mu_2((2-H(\phi))\nabla u_{h,2}&-\mathbf{g}_{h,2})\cdot\nabla w_{h,2}\xdif \bx
		+\int_{\Omega} \mathcal{E}_{\mathrm{cp}} q_2(u_{h,1},u_{h,2},w_{h,2})
		\delta_\epsilon(\phi)|\nabla\phi|\xdif \bx \notag\\ =& \int_{\Omega_1}f_1w_{h,1}\xdif \bx+ \int_{\Omega_2}f_2w_{h,2}\xdif \bx
		\quad \forall \{w_{h,1},\,w_{h,2}\}\in V_{h,1}(\delta)\times V_{h,2}(\delta).
	\end{split}
\end{equation}
In the numerical examples of Section \ref{sec:num}, we define $\delta_\epsilon(\phi)$
using the Gaussian regularization
\begin{equation}\label{eq:gauss}
	\delta_\epsilon(\phi) =H_\epsilon'(x)=
	\frac{1}{\epsilon}\sqrt{\frac{\pi}{9}}\exp\left(\frac{-\pi^2\phi^2}{9\epsilon^2}\right),\qquad
	H_\epsilon(\phi) = \frac{1}{2}\left(1+\erf\left(\frac{\pi \phi}{3\epsilon}\right)\right),
\end{equation}
where 
$\erf(x)=\frac{2}{\sqrt{\pi}}\int_0^x e^{-y^2}\mathrm{d}y$ is the Gauss error function and $H_\epsilon(\phi)$ is a smooth approximation to the sharp Heaviside function~$H(\phi)$.

\begin{rmrk}
	For $\delta=\mathrm{diam}(\bar\Omega)$, which is a valid choice, the extended subdomains $\Omega_{h,k}(\delta),\ k=1,2$
	coincide with $\Omega$, and the spaces $V_{h,k}(\delta),\ k=1,2$ coincide with the FE space for the whole domain. This version of the diffuse-interface PG method is particularly easy to implement in an existing finite element code because it boils down to solving two coupled subproblems on the same mesh. 
\end{rmrk} 

\section{Analysis of the PG stabilization}
\label{sec:analysis}

In this section, we analyze the sharp interface formulation \eqref{wsiab} of our unfitted finite
element method and prove several useful properties of the \rev{projected gradient} stabilization. \rev{The analysis covers only the lowest order case of $p=1$.} For the purpose of analysis, we assume that the mesh is shape regular and quasiuniform. Furthermore, we assume that $\delta\le C\, h$, where $\delta $ is the extension width and $C\ge0$ is an $O(1)$ constant. For brevity, we drop  $\delta$ in the notation for the extended domains and finite element spaces below.

We define the  operator  $\mathcal{P}_k:(L^2(\Omega_{h,k}))^d\to 
\left(V_{h,k}\right)^d$ for vector fields $\mathbf{v}\in(L^2(\Omega_{h,k}))^d$
as follows:
\begin{align}\label{lumproj}
	\mathcal{P}_k\mathbf{v}=\sum_{j=1}^{N_{h,k}}\mathbf{g}_{j,k}(\mathbf{v})\varphi_j,\qquad
	\mathbf{g}_{j,k}(\mathbf{v})=\frac{\int_{\Omega_{h,k}} \mathbf{v}\varphi_j\xdif \bx}{\int_{\Omega_{h,k}}\varphi_j\xdif \bx},\qquad k=1,2.
\end{align}
The bilinear form \eqref{neumannpen} of the PG stabilization term uses
$\mathbf{g}_{h,k}=\mathcal{P}_k\nabla u_{k}$ and can be written as
\begin{equation}\label{form}
	s_{h}(u,v)=\sum_{k=1}^2s_{h,k}(u_k,v_h),\quad	s_{h,k}(u,v)=\int_{\Omega_{h,k}}(\nabla u-\mathcal{P}_k\nabla u)\cdot\nabla v
\xdif \bx,\quad u,v\in H^1(\Omega_{h,k}).
\end{equation}

\rev{
  \begin{remark}
A very similar ghost penalty stabilization term was used by Badia et al.
\cite[Eq. (23)]{badia}. It differs from  \eqref{form}
in that the projection operator $\mathcal{P}_k$ is replaced by an
extension operator $\mathcal{P}_k^{ag}$ and the symmetric form
$$
s_{h,k}(u,v)=\int_{\Omega_{h,k}}\gamma(\nabla u-\mathcal{P}_k^{ag}\nabla u)\cdot
(\nabla v-\mathcal{P}_k^{ag}\nabla v)
\xdif \bx,\quad u,v\in H^1(\Omega_{h,k})
$$
of the stabilization term with a penalty parameter $\gamma$
is employed. As we show in Lemma~\ref{lemma1} 
below, our PG stabilization is also symmetric and positive semidefinite
for $p=1$.
  \end{remark}
}  

We shall write $a\lesssim b$ if $a\le c\,b$ holds with a constant $c$ that may depend only on the shape regularity of the mesh but not  on discretization parameter $h$ or the position of the interface $\Gamma$ in the background mesh. We write $a\simeq b$ if both   $a\lesssim b$ and  $b\lesssim a$ hold.
We first prove the following result.

\begin{lmm}\label{lemma1} The bilinear form  $s_{h,k}(u,v)$
is symmetric. Furthermore, the following equivalence holds 
\begin{equation}\label{form2}
	s_{h,k}(u,u)\simeq\|\nabla u-\mathcal{P}_k(\nabla u)\|^2_{ L^2(\Omega_{h,k})}+  h^2\|\nabla\mathcal{P}_k(\nabla u)\|^2_{ L^2(\Omega_{h,k})}\quad \forall\,u\in H^1(\Omega_{h,k}).
\end{equation}
\end{lmm}  
\begin{proof} Let $M=(M_{ij})_{i,j=1}^N$ be the finite element mass matrix (Gram matrix) corresponding to the nodal basis of 
the space $V_{h,k}$. Denote by $\widetilde M$ the lumped mass matrix, i.e., a diagonal matrix  with
$\widetilde M_{ii}=  \sum_{j=1}^N M_{ij}$. Note that $ M_{ij}\ge0$ for linear Lagrange finite elements and that
\begin{equation}\label{eq:aux608}
	\rw^T(\widetilde M - M)\rw=\sum_{i=1}^N\sum_{j=i+1}^NM_{ij}(\rw_i-\rw_j)^2\qquad \forall \rw\in \xR^N
\end{equation}    
holds by definition of $M$ and $\widetilde M$. It follows that
the matrix $\widetilde M - M $ is positive semi-definite.

For $w,v\in V_{h,k}$, we denote by $\rw,\rv \in \xR^N$ the corresponding vectors of coefficients in the Lagrange basis of $V_{h,k}$.  The  result below follows from \eqref{eq:aux608} and is well known (see, e.g., \cite[Proposition 1]{becker3}): 
\begin{equation}\label{eq:aux672}
	c_0\sum_{K\in\mathcal{T}_{h,k}} h_K^2\|\nabla w\|^2_{L^2(K)}\le \rw^T(\widetilde M - M)\rw\le C_0\sum_{K\in\mathcal{T}_{h,k}} h_K^2\|\nabla w\|^2_{L^2(K)},
\end{equation}
where $h_K=\mathrm{diam}(K)$ and
$C_0, c_0>0$ are constants depending only on the shape regularity of the mesh.	

A `lumped' $L^2$ inner product is defined by 
\[
(w_h,v_h)_\ell : = \rw^T\widetilde M\rv.
\] 
It is easy to verify that $\mathcal{P}_k$ satisfies the identity 
\begin{equation}\label{eq:proj}
	(\phi_h,\mathcal{P}_k g)_\ell = 	(\phi_h, g)_{L^2(\Omega_{h,k})} \quad \forall\, \phi_h\in V_{h,k}, ~g\in L^2(\Omega_{h,k}).
\end{equation}
Substituting the scalar components of $\boldsymbol{\phi}_h=\mathcal{P}_k\nabla v$ and $\mathbf {g}=\nabla u$  into
\eqref{eq:proj}, we find that
\begin{equation}\label{eq:proj2}
	(\mathcal{P}_k\nabla v,\mathcal{P}_k\nabla u)_\ell = (\mathcal{P}_k\nabla v,\nabla u)_{L^2(\Omega_{h,k})}.
\end{equation}
Here and below, the
inner products should be understood componentwise. Using \eqref{eq:proj2}, we obtain
\[
\begin{split}
	s_{h,k}(u,v)&=(\nabla u-\mathcal{P}_k\nabla u,\nabla v-\mathcal{P}_k\nabla v)_{L^2(\Omega_{h,k})} +
	(\nabla u-\mathcal{P}_k\nabla u,\mathcal{P}_k\nabla v)_{L^2(\Omega_{h,k})} \\
	&=(\nabla u-\mathcal{P}_k\nabla u,\nabla v-\mathcal{P}_k\nabla v)_{L^2(\Omega_{h,k})} +
	(\mathcal{P}_k\nabla v,\mathcal{P}_k\nabla u)_\ell-(\mathcal{P}_k\nabla u,\mathcal{P}_k\nabla v)_{L^2(\Omega_{h,k})}.
\end{split}
\]
This shows the symmetry of $s_{h,k}(u,v)$. 

To prove \eqref{form2}, we show that the difference of the last two terms on the right-hand side is equivalent to $h^2\|\nabla\mathcal{P}_k(\nabla u)\|^2_{ L^2(\Omega_{h,k})}$.  This follows from \eqref{eq:aux672}, the quasiuniformity of the mesh, and the observation that $\mathcal{P}_k\nabla u$ is an element of $(V_{h,k})^d$. Indeed, for any $w_h\in V_{h,k}$ we have the estimate
\[
(w_h,w_h)_\ell-(w_h,w_h)_{L^2(\Omega_{h,k})} = \rw^T\widetilde M\rw-\rw^T M\rw = \rw^T(\widetilde M - M)\rw \simeq \sum_{K\in\mathcal{T}_h} h_K^2\|\nabla w_h\|^2_{L^2(K)}, 
\]  	 
which holds componentwise for $\mathbf w_h=\mathcal{P}_k(\nabla u_h)$. This proves  \eqref{form2}.
\end{proof}

Let us now estimate $s_{h,k}(u,u)$ from below by a sum of local stabilization forms associated with element patches. We denote by $\mathcal{F}_{h,k}$ the set of internal faces (edges in 2D) of the sub-triangulation $\mathcal T_{h,k}$ and let $\omega(f)=K\cup K'$ be the union of two simplices sharing a face $f\in\mathcal{F}_{h,k}$. Using the notation {$h_f=\max\{h_{K},h_{K'}\}$} for the local mesh size of the patch, we define 
\[
s_{h,k}^f(u,v) = \int_{\omega(f)}[(\nabla u-\mathcal{P}_k\nabla u)\cdot(\nabla v-\mathcal{P}_k\nabla v)]\xdif \bx + h^2_f\int_{\omega(f)}(\nabla \mathcal{P}_k\nabla u):(\nabla \mathcal{P}_k\nabla v)\xdif \bx.
\]

The estimate \eqref{form2} implies
\begin{equation}\label{eq:aux698}
s_{h,k}(u,u)\gtrsim \sum_{f\in \mathcal{F}_{h,k}} s_{h,k}^f(u,u)\quad \forall\,u\in H^1(\Omega_{h,k})
\end{equation}
with some $c>0$ depending only on the shape regularity of the mesh.	

We are now ready to prove a key technical result. 

\begin{lmm} \label{L:localCoerc}
For arbitrary $K_1, K_2 \in{\mathcal T}_{h,k}$ such that 
$K_1\cup K_2= \omega(f)$ for a face $f\in{\mathcal T}_{h,k}$, the following estimates hold
\begin{eqnarray}
	c_1\|\nabla u\|_{L^2(K_2)}^2\le \| \nabla u\|_{L^2(K_1)}^2 +  s_{h,k}^f(u,u)\quad \forall\,u\in H^1(\Omega_{h,k}) \label{est1},\\
	c_0\|u_h\|_{L^2(K_2)}^2\le  \| u_h\|_{L^2(K_1)}^2 + h^{2} s_{h,k}^f(u_h,u_h)\quad \forall\,u_h\in V_{h,k} \label{est0}
\end{eqnarray}
with some constants $c_0,c_1>0$ depending only on the shape regularity of the mesh.
\end{lmm} 
\begin{proof} By the triangle inequality, we have
\begin{equation}\label{eq:aux739}
	\|\nabla u\|_{L^2(K_2)}\le \|\nabla u- \mathcal{P}_k\nabla u\|_{L^2(K_2)}+\|\mathcal{P}_k\nabla u\|_{L^2(K_2)}.
\end{equation}
Let $\mathbf w_h=\mathcal{P}_k\nabla u$. Then $\mathbf w_h\in (V_{h,k})^d$. Using the \rev{Friedrichs inequality and a scaling argument}, we estimate the second term on the right-hand side of \eqref{eq:aux739} as follows: 
\begin{equation*}
	\begin{split}
		\|\mathcal{P}_k\nabla u\|_{L^2(K_2)}^2&=\|\mathbf w_h\|_{L^2(K_2)}^2\lesssim   h^2_{K_2}\|\nabla\mathbf w_h\|_{L^2(K_2)}^2 +  h_f\|\mathbf w_h\|_{L^2(f)}^2 \\
		& \lesssim   h^2_{K_2}\|\nabla \mathbf w_h\|_{L^2(K_2)}^2 + h^2_{K_1}\|\nabla \mathbf w_h\|_{L^2(K_1)}^2 + \|\mathbf w_h\|_{L^2(K_1)}^2\\ 
		& \lesssim   h^2_f\| \nabla \mathcal{P}_k\nabla u \|^2_{\omega(f)} +  \| \nabla u\|_{L^2(K_1)}^2. 
	\end{split}
\end{equation*}
Together with \eqref{eq:aux739}  and the definition of $s_{h,k}^f(\cdot,\cdot)$ this proves the result claimed in \eqref{est1}. 

To prove the validity of \eqref{est0}, we use the  finite element  inverse and trace  estimates, as well as \eqref{est1}, in the chain of estimates
\begin{equation*}
	\begin{split}
		\|u_h\|_{L^2(K_2)}^2&\lesssim   h^2_{K_2}\|\nabla u_h\|_{L^2(K_2)}^2 +  h_f\|u_h\|_{L^2(f)}^2 \\
		&\le    h^2_{K_2}c_1^{-1}( \| \nabla u_h\|_{L^2(K_1)}^2 +  s_{h,k}^f(u_h,u_h)) +  h_f\|u_h\|_{L^2(f)}^2 \\
		&\lesssim    h^2_{K_2}s_{h,k}^f(u_h,u_h) +  h^2_{K_1}\|\nabla u_h\|_{L^2(K_1)}^2 + \|u_h\|_{L^2(K_1)}^2 \\
		&\lesssim  h^2_{K_2}s_{h,k}^f(u_h,u_h) + \|u_h\|_{L^2(K_1)}^2,		
	\end{split}
\end{equation*}
which complete the proof of the lemma.
\end{proof}

Recall that $\Omega_{h,k}$ consists of simplices contained in $\bar\Omega_k$ and layers of
simplices $K\not\subset\bar\Omega_k$ that intersect the $\delta$-neighborhood of $\Gamma$. Because of our assumption on $\delta$, the number
of layers is finite and independent of $h$.
Then the smoothness of $\partial\Omega$ and the regularity of the mesh imply (see, e.g., \cite{olsh}) that any $K\in\mathcal{T}_h\cap\Omega_{h,k}$ can be reached from a strictly interior simplex by crossing a `finite' number of faces. 
Moreover, the number of paths originating from any interior simplex is uniformly bounded. Using this argument, as well as the estimates \eqref{est1} and \eqref{eq:aux698}, we find that
\begin{equation}\label{coer}
\|\nabla u\|^2_{L^2(\Omega_{h,k})} \lesssim \|\nabla u\|^2_{L^2(\Omega_k)} + s_{h,k}(u,u)\quad\forall u\in H^1(\Omega_{h,k}).
\end{equation}

This coercivity result is  important for performing stability and error analysis.  The estimate provided by \eqref{est0} is helpful for the analysis of  evolving interface problems, which is not a subject of this paper.
\smallskip

We now proceed with the consistency analysis.  Consistency results for the elliptic part and the Nitsche terms in the finite element formulation are provided by the standard theory.  Therefore, we focus our attention on the new stabilization--extension term. 

We first need several results regarding the stability and approximation properties of the operator $\mathcal{P}_k$
defined in \eqref{lumproj}. 

\begin{lmm}\label{ProjProp}
The mapping $\mathcal{P}_k$ is both  $L^2$- and $H^1$-stable, i.e., it has the property that
\begin{equation}\label{Pstab}
	\|\mathcal{P}_k v\|_{L^2(\Omega_{h,k})}\lesssim 	\|v\|_{L^2(\Omega_{h,k})}\quad \text{and}\quad 
	\|\nabla \mathcal{P}_k w\|_{L^2(\Omega_{h,k})}\lesssim 	\|\nabla w\|_{L^2(\Omega_{h,k})}
\end{equation} 
for any $ v\in L^2(\Omega_{h,k})$, $ w\in H^1(\Omega_{h,k})$, $k=1,2$. Furthermore, it holds
\begin{equation}\label{Papprox}
	\|v- \mathcal{P}_k v\|_{L^2(\Omega_{h,k})}\lesssim  h\|\nabla v\|_{L^2(\Omega_{h,k})}\quad \forall\, v\in H^1(\Omega_{h,k}).
\end{equation} 
\end{lmm}
\begin{proof}
The consistent and lumped mass matrices are spectrally equivalent~\cite{Wathen}, which implies $\|\phi_h\|_\ell\lesssim \|\phi_h\|_{L^2(\Omega_{h,k})}$ for any $\phi_h\in V_{h,k}$.  Using this fact and choosing $\phi_h=\mathcal{P}_k v$ in \eqref{eq:proj}, we obtain
\begin{equation*}
	\|\mathcal{P}_k v\|_{L^2(\Omega_{h,k})}^2
	\lesssim 	\|\mathcal{P}_k v\|_{\ell}^2= 	(\mathcal{P}_k v, v)_{L^2(\Omega_{h,k})} \le 
	\|\mathcal{P}_k v\|_{L^2(\Omega_{h,k})}	\|v\|_{L^2(\Omega_{h,k})}.
\end{equation*}
This implies the first inequality in \eqref{Pstab}.

In the next step, we first check the validity of \eqref{Papprox} assuming that $v_h\in V_{k,h}\subset H^1(\Omega_{h,k})$.
In this case, due to the  assumption of a quasiuniform mesh, we have 
\begin{equation}\label{eq:aux726}
	\begin{split}
		\|v_h - \mathcal{P}_k v_h\|_{L^2(\Omega_{h,k})}^2 &
		\simeq h^d \sum_{j=1}^{N_{h,k}} \left[\frac{\int_{\omega(\bx_j)} v_h\varphi_j\xdif \bx}{\int_{\omega(\bx_j)}\varphi_j\xdif \bx} - v_h(\bx_j) \right]^2  \\
		&\simeq h^{-d} \sum_{j=1}^{N_{h,k}} \left[ \int_{\omega(\bx_j)}v_h\varphi_j\xdif \bx  - v_h(\bx_j) \int_{\omega(\bx_j)}\varphi_j\xdif \bx \right]^2.
	\end{split}
\end{equation}
Here $\omega(\bx_j)$ denotes a domain consisting of simplices that share the grid node $\bx_j$.
By a standard argument that involves a pullback to a reference triangle, using a norm equivalence  and  a pushforward, one can show that for any $K\in\mathcal{T}_h$ and a grid node $\bx_j\in K$  the following estimate holds: 
\[
\Big|\int_{K} v_h\varphi_j\xdif \bx - v_h(\bx_j)\int_{K}\varphi_j \xdif \bx\Big|\lesssim 
h \|\nabla v_h\|_{L^2(K)}\|\varphi_j\|_{L^2(K)}.
\]
Summing up over all simplices from $\omega(\bx_j)$ and   noting that $\|\varphi_j\|_{L^2(\omega(\bx_j))} \simeq h^{d/2}$ gives
\[
\Big|\int_{\omega(\bx_j)} v_h\varphi_j\xdif \bx - v_h(\bx_j)\int_{\omega(\bx_j)}\varphi_j \xdif \bx\Big|\lesssim 
h^{d/2+1} \|\nabla v_h\|_{L^2(\omega(\bx_j))}.
\]
Substituting this estimate into the right hand side of \eqref{eq:aux726}  proves  \eqref{Papprox} for $v_h\in V_{k,h}$.
To extend \eqref{Papprox}  to arbitrary $v \in H^1(\Omega_{h,k})$, we consider the $L^2$-orthogonal projection $\mathcal Q_h v \in V_{k,h} $ and apply the triangle inequality to obtain an upper bound for
\[
\|v- \mathcal{P}_k v\|_{L^2(\Omega_{h,k})} \le 
\|\mathcal{P}_k (\mathcal Q_h v - v)\|_{L^2(\Omega_{h,k})}+
\|\mathcal{P}_k (\mathcal Q_h v) - \mathcal Q_h  v\|_{L^2(\Omega_{h,k})} +
\|v- \mathcal Q_h v\|_{L^2(\Omega_{h,k})}.
\]
We then use the $L^2$-stability of $\mathcal{P}_k$, estimate \eqref{Papprox} on $V_{k,h}$, as well as the $H^1$-stability and approximation properties of the $L^2$-orthogonal projection on quasiuniform meshes to show that the right hand side is bounded by $h\|\nabla v\|_{L^2(\Omega_{h,k})}$, which proves that \eqref{Papprox} is valid for all $v \in H^1(\Omega_{h,k})$. 

The second estimate in \eqref{Pstab} can now be deduced using the triangle inequality, the finite element inverse inequality, \eqref{Papprox}, as well as the $H^1$-stability and approximation properties of the $L^2$-orthogonal projection:
\[
\begin{split}
	\|\nabla\mathcal{P}_k v\|_{L^2(\Omega_{h,k})} & \le 
	\|\nabla (\mathcal{P}_k v - \mathcal Q_h  v)\|_{L^2(\Omega_{h,k})} +
	\|\nabla \mathcal Q_h v\|_{L^2(\Omega_{h,k})} \\
	&\lesssim
	h^{-1}\|\mathcal{P}_k v - \mathcal Q_h  v\|_{L^2(\Omega_{h,k})} +
	\|\nabla v\|_{L^2(\Omega_{h,k})} \\
	& \lesssim
	h^{-1}\|\mathcal{P}_k v - v \|_{L^2(\Omega_{h,k})}+  h^{-1}\|v - \mathcal Q_h  v\|_{L^2(\Omega_{h,k})} +
	\|\nabla v\|_{L^2(\Omega_{h,k})} \\
	& \lesssim 
	\|\nabla v\|_{L^2(\Omega_{h,k})}.
\end{split}
\]

\end{proof}

Next, we prove the following continuity and consistency result for   the PG stabilization operator. 

\begin{lmm}\label{lemma4} There holds
\begin{eqnarray}\label{sc1}
	{s}_{h,k}(u,v)&\lesssim&  \|u\|_{H^1(\Omega_{h,k})}\|v\|_{H^1(\Omega_{h,k})}\quad \forall\,u,v\in H^1(\Omega_{h,k}),\, v\in H^1(\Omega_{h,k}),\\
	{s}_{h,k}(u,u)&\lesssim& h^2 \|u\|_{H^2(\Omega_{h,k})}^2\quad \forall\,u\in H^{2}(\Omega_{h,k}).
	\label{sc2}
\end{eqnarray}
\end{lmm}  
\begin{proof} The continuity estimate \eqref{sc1} follows from the Cauchy--Schwarz inequality and the $L^2$-stability of $ {\mathcal{P}}_k$: 
\begin{equation*}
	\begin{split}
		{s}_{h,k}(u,v)&=\int_{\Omega_{h,k}}(\nabla u - {\mathcal{P}}_k\nabla u)\cdot\nabla v\,\xdif \bx \le   \|\nabla u - {\mathcal{P}}_k\nabla u\|_{L^2(\Omega_{h,k})}\|\nabla v\|_{L^2(\Omega_{h,k})}\\
		&\le 
		(\|\nabla u\|_{L^2(\Omega_{h,k})}+\| {\mathcal{P}}_k\nabla u\|_{L^2(\Omega_{h,k})})\|\nabla v\|_{L^2(\Omega_{h,k})}
		\lesssim  \|u\|_{H^1(\Omega_{h,k})}\|v\|_{H^1(\Omega_{h,k})}.
	\end{split}
\end{equation*}
To verify \eqref{sc2}, we employ \eqref{form2},  the $H^1$ stability of  ${\mathcal{P}}_k$ and \eqref{Papprox}:
\begin{equation*}
	s_{h,k}(u,u)\lesssim \|\nabla u-\mathcal{P}_k(\nabla u)\|^2_{ L^2(\Omega_{h,k})}+ c\, h^2\|\nabla\mathcal{P}_k(\nabla u)\|^2_{ L^2(\Omega_{h,k})} \lesssim  h^2\|\nabla^2 u\|^2_{ L^2(\Omega_{h,k})}.
\end{equation*}
\end{proof}

To formulate a convergence result, we need the following norm 
\[
\|v\|_\ast^2 = \sum_{k=1}^2 \left\{ 
\| v_k\|_{H^1(\Omega_{h,k})}^2
+ h^{-1}\|v_k\|_{L^2(\Gamma)}^2+ {s}_{h,k}(v_k,v_k)
\right\}\quad \text{for}~ v\in H^1(\Omega_{h,1})\times H^1(\Omega_{h,2}).
\]

If $u_k\in H^1(\Omega_k)$, there exists a bounded linear operator $\mathcal{E}=\mathcal{E}_k:H^1(\Omega_k)\to H^1(\mathbb{R}^d)$ that extends $u_k$ to $\mathcal{E}u_k\in H^1(\Omega)$. Moreover, $\mathcal{E}_k$ is also bounded as an operator from $H^2(\Omega_k)$ to $H^2(\mathbb{R}^d)$; see~\cite{Stein}. Whenever this causes no confusion, we will identify $u_k$ with its extension $\mathcal{E}u_k$ to $\Omega$. 

 The theorem below proves the main convergence result.

\begin{thrm}\label{th1} 
	\rev{Assume that $u\in H^2(\Omega_{1})\times  H^2(\Omega_{2})$ solves problem \eqref{eq:inter}. 
	We associate with $u$ its double-extension $\mathcal{E}u=\{\mathcal{E}u_1,\mathcal{E}u_2\}$. Assume   $\delta\lesssim h$ for the extension width $\delta$.} Let $u_h$ be a solution to \eqref{wsiab}. Then
\rev{
\begin{align}\label{err_est1}
	\|u-u_h\|_\ast &\lesssim 
	h\sum_{k=1}^2\|u\|_{H^2(\Omega_{k})}, 
\\	\label{err_est2}
\|u-u_h\|_{L^2(\Omega)} &\lesssim 
h^2\sum_{k=1}^2\|u\|_{H^2(\Omega_{k})}. 
\end{align}
}
\end{thrm}
\begin{proof}
Let $I_k$ denote the nodal interpolation operator for $\Omega_{h,k}$ and
$I_h(u)\in V_{h,1}\times V_{h,2}$ the componentwise nodal interpolant of $u$, i.e., $I_h(u)=\{I_1(\mathcal{E}u_1),I_2(\mathcal{E}u_2)\}$. 
Owing to \eqref{coer}, one can follow standard arguments (see, e.g., \cite{hansbo2,burman2012,cutfem}) to check that the bilinear form
$$A_h(u,v):=a(u,v)+s_h(u,v)$$
is uniformly coercive and $a$ is continuous on the finite element space with respect to the $\|\cdot\|_\ast$ norm, \rev{i.e. $A_h(v_h,v_h)\gtrsim \|v_h\|^2_\ast$, for all $v_h\in  V_{h,1}\times V_{h,2}$ and $a(u_h,v_h)\lesssim \|u_h\|_\ast\|v_h\|_\ast$ for all $u_h, v_h\in  V_{h,1}\times V_{h,2}$}. We therefore have 
\begin{equation*}
	\begin{split}
		\|I_h(u)-u_h\|_\ast^2&\lesssim A_h(I_h(u)-u_h,I_h(u)-u_h)\\ 
		&=  A_h(I_h(u)-u,I_h(u)-u_h) + s_h(u,I_h(u)-u_h) \\ 
		&= a(I_h(u)-u,I_h(u)-u_h)+ s_h(I_h(u)-u,I_h(u)-u_h)+ s_h(u,I_h(u)-u_h).
	\end{split}
\end{equation*}
\rev{Note that due to symmetry and \eqref{form2}, the bilinear form $s_h(\cdot,\cdot)$ defines a semi-inner product. We use this property to apply the Cauchy--Schwarz inequality to the last term on the right hand side.} 
Using \eqref{sc1}, \rev{\eqref{sc2}},  \rev{the} continuity estimate for the bilinear form  $a(\cdot,\cdot)$, the interpolation error estimate \rev{(see, e.g., \cite[Lemma~5]{burman2012}),
\[
	\|u-I_h(u)\|_\ast \lesssim h \sum_{k=1}^2\|\mathcal{E}u\|_{H^2(\Omega_{h,k})},
\] 
}
and the $H^2$ boundedness of the extension operator $\mathcal{E}$, we find 
\begin{equation}\label{aux767}
	\|I_h(u)-u_h\|_\ast \lesssim h \sum_{k=1}^2\|\mathcal{E}u\|_{H^2(\Omega_{h,k})} \lesssim h \sum_{k=1}^2\|u\|_{H^2(\Omega_{k})}.
\end{equation}
This estimate, the triangle inequality, and interpolation properties prove \eqref{err_est1}.

We now proceed with a duality argument and consider the solution to the problem 
\begin{equation}\label{eq:dual}
	\begin{aligned}
		-\nabla\cdot(\mu\nabla z)&=u-u_h\quad\text{in }\Omega,\\
		z &= 0\quad\text{on }\partial\Omega,\\
		[\![ z]\!]&=0\quad\text{on }\Gamma,\\
		[\![ \mu\nabla z]\!]&=0\quad\text{on }\Gamma.
	\end{aligned}
\end{equation}
It is well known  \cite{zou} that $z\in  H^2(\Omega_{1})\times  H^2(\Omega_{2})$ and 
\begin{equation}\label{eq:aux880}
	\sum_{k=1}^2\|z\|_{H^2(\Omega_{k})} \lesssim \|u-u_h\|_{L^2(\Omega)}. 
\end{equation}
Denote by $I_h(z)\in V_{h,1}\times V_{h,2}$ the componentwise nodal interpolant for (the extension of) $z$. Testing the first equation in \eqref{eq:dual} by $u-u_h$, integrating by parts, and utilizing \eqref{wsiab}  one finds the relation
\[
\|u-u_h\|_{L^2(\Omega)}^2= a(u-u_h,z-I_h(z))+ s_h(u_h,I_h(z)). 
\] 
\rev{Recall that the bilinear form $s_h(\cdot,\cdot)$ defines a semi-inner product.} We use this property below to apply the Cauchy--Schwarz and triangle inequalities. 
Exploiting the continuity, the interpolation properties of $I_h(z)$ and \eqref{eq:aux880}, we arrive at  
\begin{equation*}
	\begin{split}
		\|u-u_h\|_{L^2(\Omega)}^2&\lesssim \|u-u_h\|_\ast\|z-I_h(z)\|_\ast + s_h^{\frac12}(u_h,u_h)s_h^{\frac12}(I_h(z),I_h(z))\\
		&\lesssim h \|u-u_h\|_\ast\|u-u_h\|_{L^2(\Omega)} + s_h^{\frac12}(u_h,u_h)s_h^{\frac12}(I_h(z),I_h(z)) . 
	\end{split}
\end{equation*}
This result and the estimate \rev{\eqref{err_est1} imply} that 
\begin{equation} \label{eq1}
	\begin{split}
		\|u-u_h\|_{L^2(\Omega)}^2
		&\lesssim\Big( h^2\sum_{k=1}^2\|u\|_{H^2(\Omega_{k})}\Big) \|u-u_h\|_{L^2(\Omega)} + s_h^{\frac12}(u_h,u_h)s_h^{\frac12}(I_h(z),I_h(z)) . 
	\end{split}
\end{equation}
We estimate the last term using the triangle  inequality,  \eqref{sc1} and \eqref{sc2} to deduce
\begin{equation}
	s_h^{\frac12}(u_h,u_h)\rev{\lesssim s_h^{\frac12}(u-u_h,u-u_h)+s_h^{\frac12}(u,u)}
	\lesssim \|u-u_h\|_\ast + h\sum_{k=1}^2\|u\|_{H^2(\Omega_{h,k})}
\end{equation}
and
\begin{equation}\label{eq3}
	\begin{split}
		s_h^{\frac12}(I_h(z),I_h(z))&\le s_h^{\frac12}(I_h(z)-z,I_h(z)-z)+s_h^{\frac12}(z,z)\\
		&\lesssim \sum_{k=1}^2(\|I_h(z)-z\|_{H^1(\Omega_{h,k})}+ h\|z\|_{H^2(\Omega_{h,k})})\\
		&\lesssim h\sum_{k=1}^2\|z\|_{H^2(\Omega_{h,k})}
		\lesssim h\sum_{k=1}^2\|z\|_{H^2(\Omega_{k})}\\
		&\lesssim  h \|u-u_h\|_{L^2(\Omega)},
	\end{split}
\end{equation}
where the last two estimates follow from the boundedness of the extension operator $\mathcal E$ in $H^2$ and from the stability result \eqref{eq:aux880}, respectively.

Collecting the estimates \eqref{eq1}--\eqref{eq3} \rev{ together with \eqref{err_est1} proves  \eqref{err_est2}}.
\end{proof}

\rev{
  \begin{remark}[Some algebraic properties]\label{rem41}
   Adopting the notation of Remark~\ref{rem22} and using a standard nodal basis again, we define $A$ as the matrix associated with the (non-stabilized) bilinear form $a(\cdot,\cdot)$.
%
Then the matrix of the finite element problem for the sharp interface formulation~\eqref{wsiab} can be written as
\[
A + L - B^T \widetilde{M}^{-1} B.
\]
Since $B^T \widetilde{M}^{-1} B$ is positive semi-definite, it follows that
$
A + L - B^T \widetilde{M}^{-1} B \le A + L
$
in the spectral sense. From the coercivity of the bilinear form $a(\cdot,\cdot)$ and the definition of the $\|\cdot\|_*$ norm, it follows that there exists a constant $c > 0$, independent of $h$ and of the position of the interface within the mesh, such that the following spectral equivalence holds:
\begin{equation}
	\label{eqSp}
	c (A + L) \le A + L - B^T \widetilde{M}^{-1} B \le A + L.
\end{equation}

In particular, this implies that on a quasi-uniform mesh, the condition number of $A + L - B^T \widetilde{M}^{-1} B$ scales like $O(h^{-2})$. Moreover, the matrix $A + L$ has a standard sparsity pattern, which can be effectively exploited in constructing suitable preconditioners for the matrix of the finite element problem. 
\end{remark}
}

\section{Test problems and numerical results}
\label{sec:num}
In this section, we perform numerical studies of our unfitted FEM with global \rev{projected gradient} stabilization for 
the elliptic interface problem \eqref{eq:inter} in two space dimensions. In our numerical experiments,
we consider different choices of diffusion coefficients $\mu$, right hand sides $f$, and boundary conditions.
The $L^2$ error of an unfitted finite element approximation in a domain $\Omega$ is given by 
\begin{align}
	\|u_h-u\|_{0,\Omega}=\sqrt{\int_{\Omega}(H(\phi)u_{h,1}+(1-H(\phi))u_{h,2}-u)^2\xdif \bx}.
\end{align}

Computations for all test cases are performed using a  C\texttt{++} implementation of the methods under investigation in the open-source finite element library MFEM \cite{mfem, mfem2}. The numerical solutions are visualized using the open-source software {\sc GlVis} \cite{glvis}.

\subsection{Poisson problem with a straight interface and a smooth solution}
\label{sec:test1dsmooth}

To begin, we apply the proposed methods to a quasi-1D version of problem \eqref{eq:inter} with a straight interface and a smooth solution. Let $\Omega=(0,1)^2$ and $\Gamma=\{(x,y)\in\bar\Omega\,:\,x=0.51\}$. The two subdomains are given by
$\Omega_1=\{(x,y)\in\Omega:x<0.51\}$ and $\Omega_2=\{(x,y)\in\Omega:x>0.51\}$. 
We set $\mu_1=10^{-8}$, $\mu_2=1$ and $f_1=2\cdot 10^{-8}$, $f_2=2$ in this test. The Dirichlet boundary conditions imposed at points $(x,y)\in\partial\Omega$ with $x=0$ or $x=1$ correspond to the analytical solution
\begin{equation*}
	u_1(x,y)=(x-0.01)(1.01-x)=u_2(x,y).
\end{equation*}
At points $(x,y)\in\partial\Omega$ with $y=0$ or $y=1$, we impose homogeneous Neumann boundary conditions.

\begin{figure}[h!]
	\centering\vspace{0.25cm}
	\begin{minipage}[t]{0.5\textwidth}\centering
		\frame{\includegraphics[angle=90,origin=c,width=0.65\textwidth]{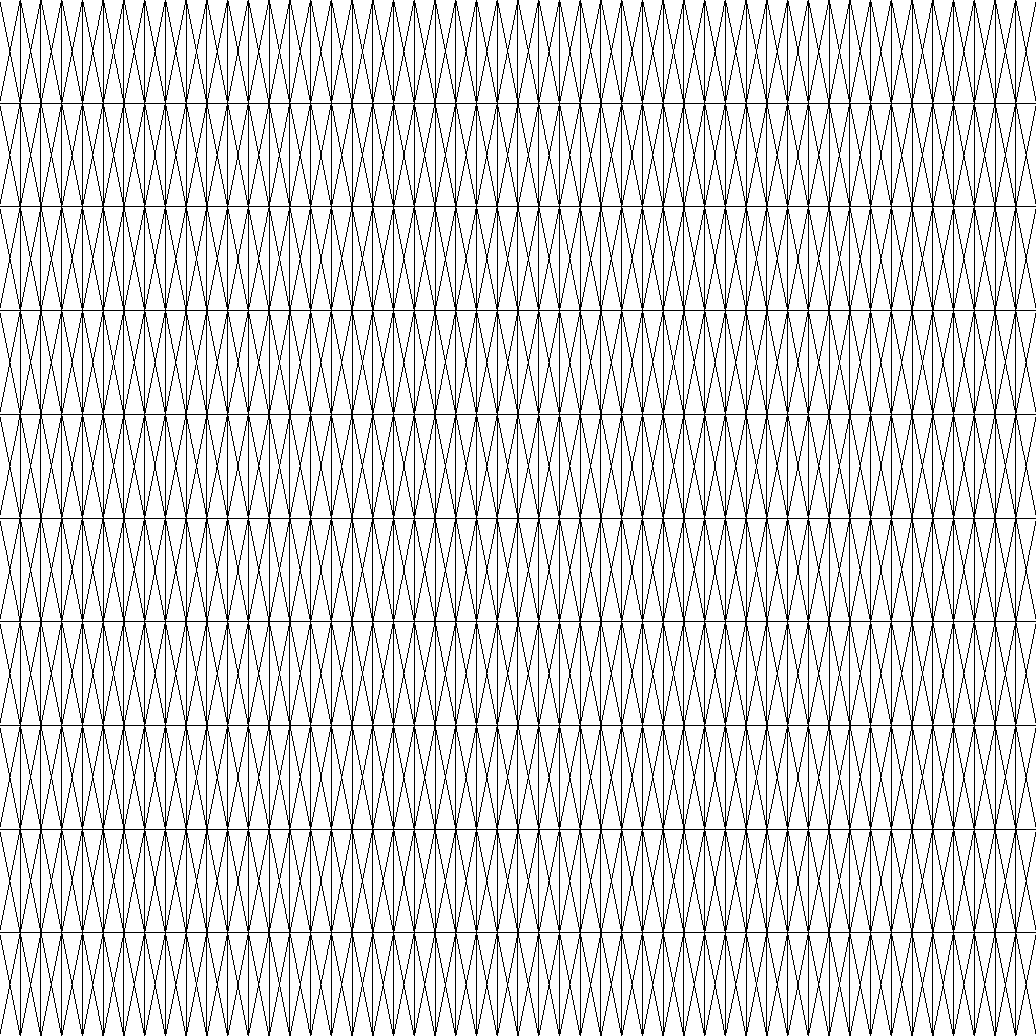}}
	\end{minipage}%
	\vskip0.25cm 
	\caption{Coarse quasi-uniform mesh with $\Delta x = 0.1$ and $\Delta y = 0.02$.}
	\label{fig:coarse}
\end{figure}

\begin{figure}[h!]
	\centering
	\begin{minipage}[t]{0.3\textwidth}
		\centering (a)  extended $u_1$\vskip0.5cm
		
		\includegraphics[width=0.9\textwidth]{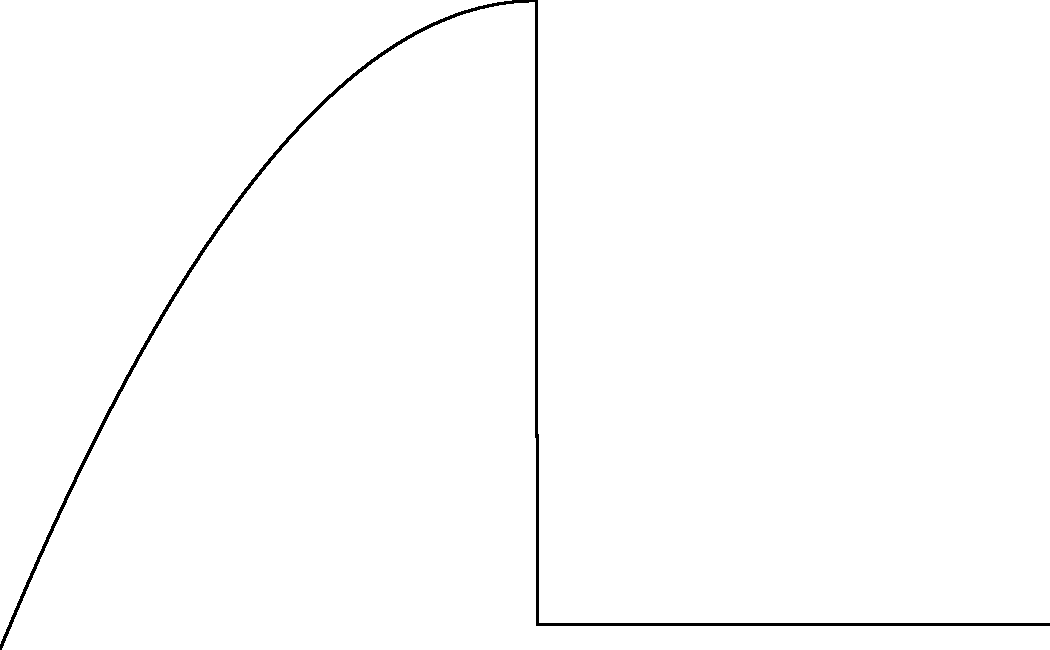}
	\end{minipage}%
	\begin{minipage}[t]{0.3\textwidth}
		\centering (b)  $u$\vskip0.5cm
		
		\includegraphics[width=0.9\textwidth]{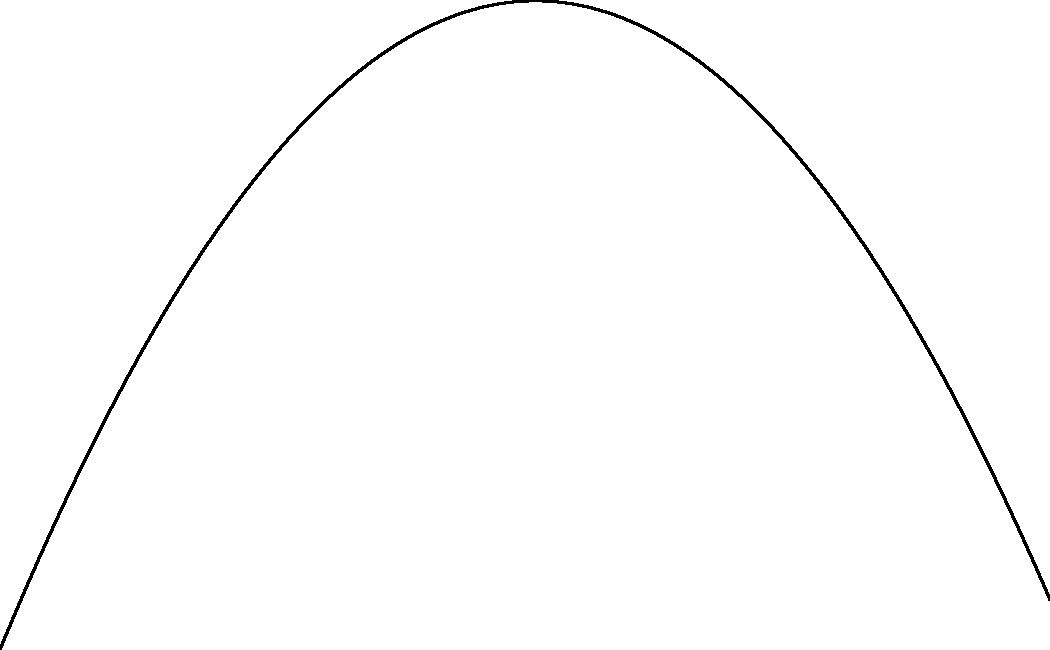}
	\end{minipage}
	\begin{minipage}[t]{0.3\textwidth}
		\centering (c)  extended $u_2$\vskip0.5cm
		
		\includegraphics[width=0.9\textwidth]{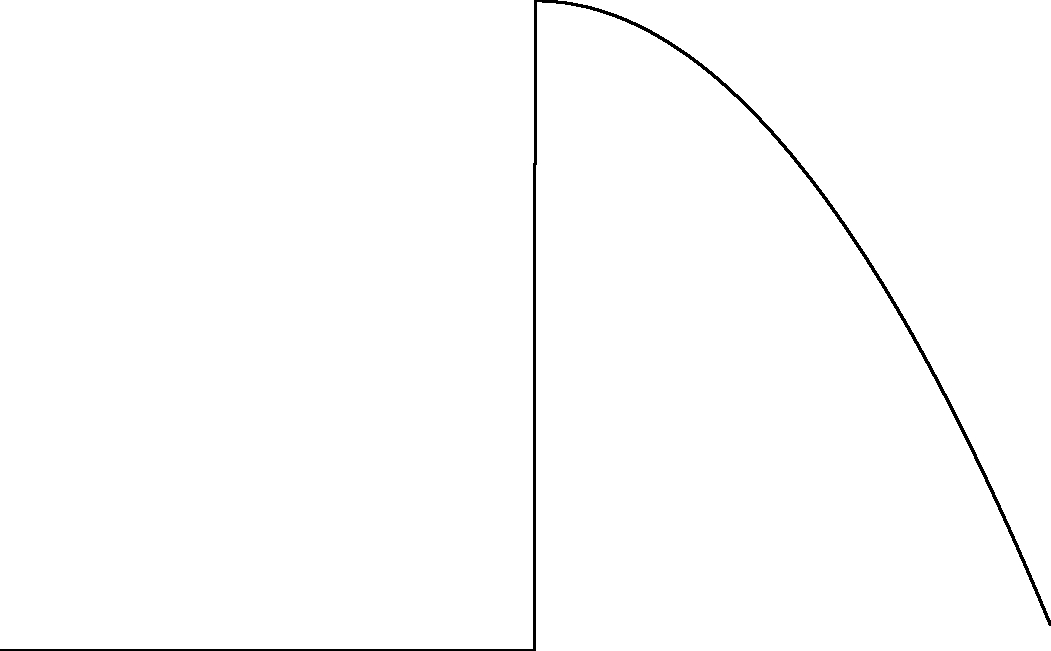}
	\end{minipage}
	\vskip0.5cm
	
	\begin{minipage}[t]{0.3\textwidth}
		\centering
		
		\includegraphics[width=0.9\textwidth]{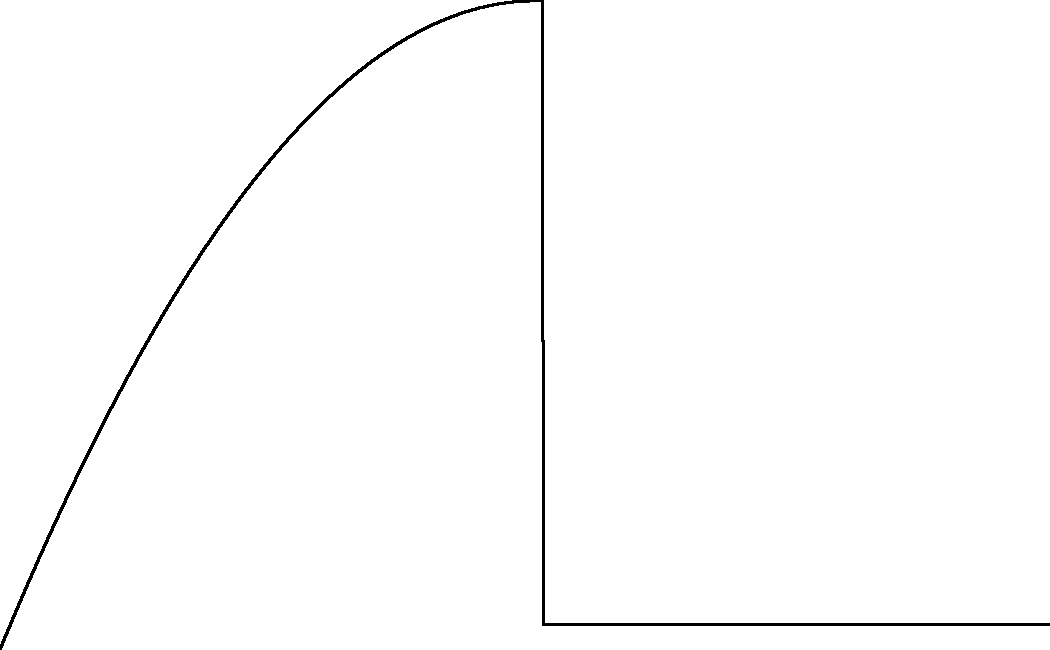}
	\end{minipage}%
	\begin{minipage}[t]{0.3\textwidth}
		\centering
		
		\includegraphics[width=0.9\textwidth]{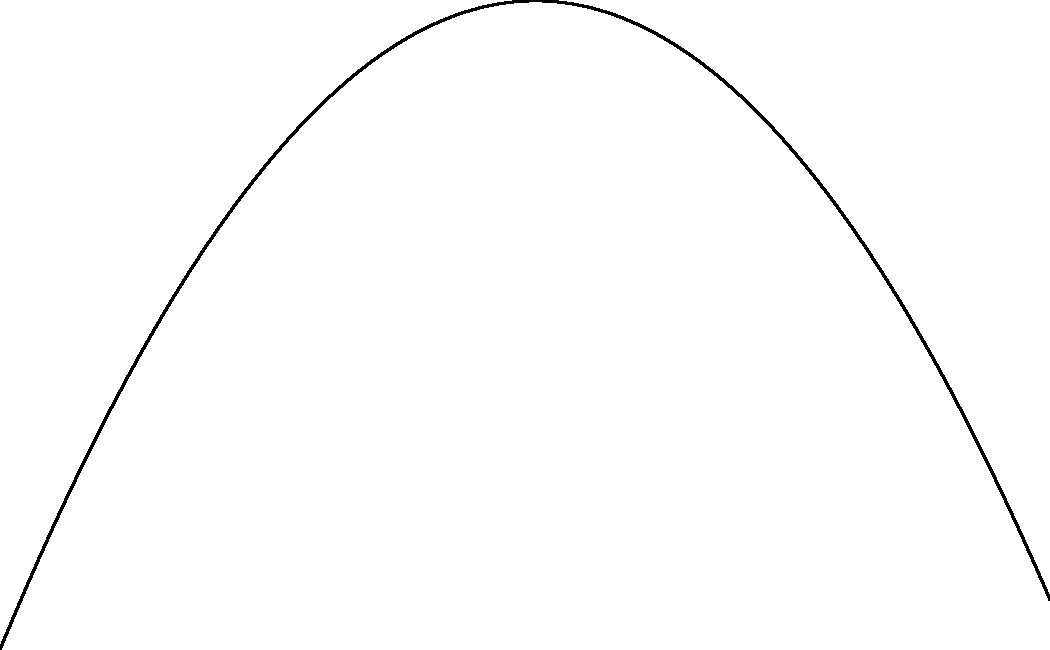}
	\end{minipage}
	\begin{minipage}[t]{0.3\textwidth}
		\centering
		
		\includegraphics[width=0.9\textwidth]{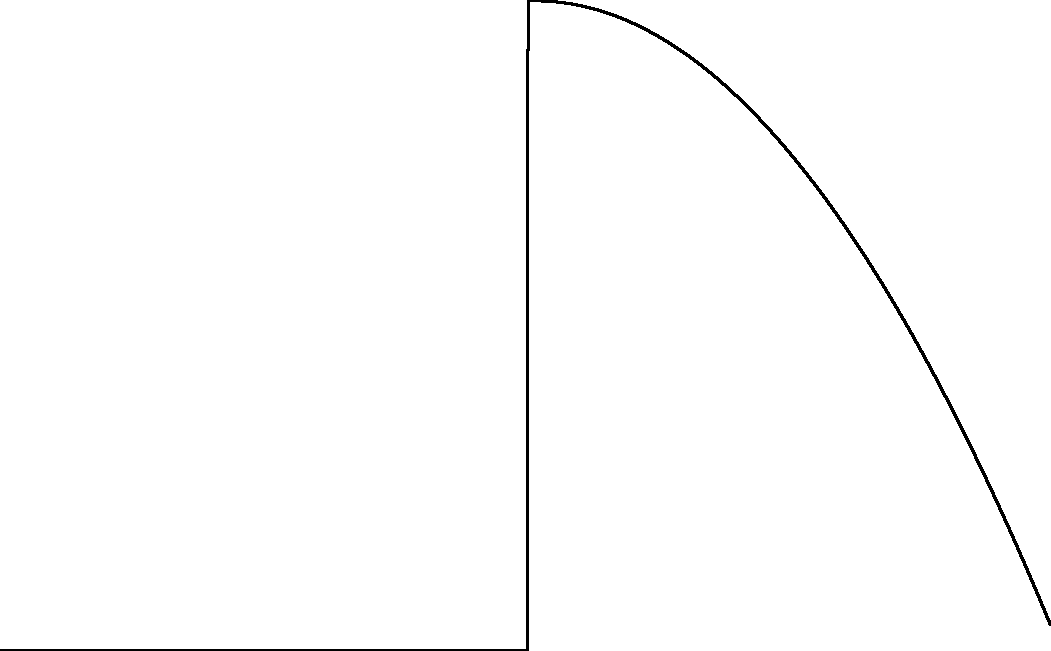}
	\end{minipage}
	\vskip0.5cm
	
	\begin{minipage}[t]{0.3\textwidth}
		\centering
		
		\includegraphics[width=0.9\textwidth]{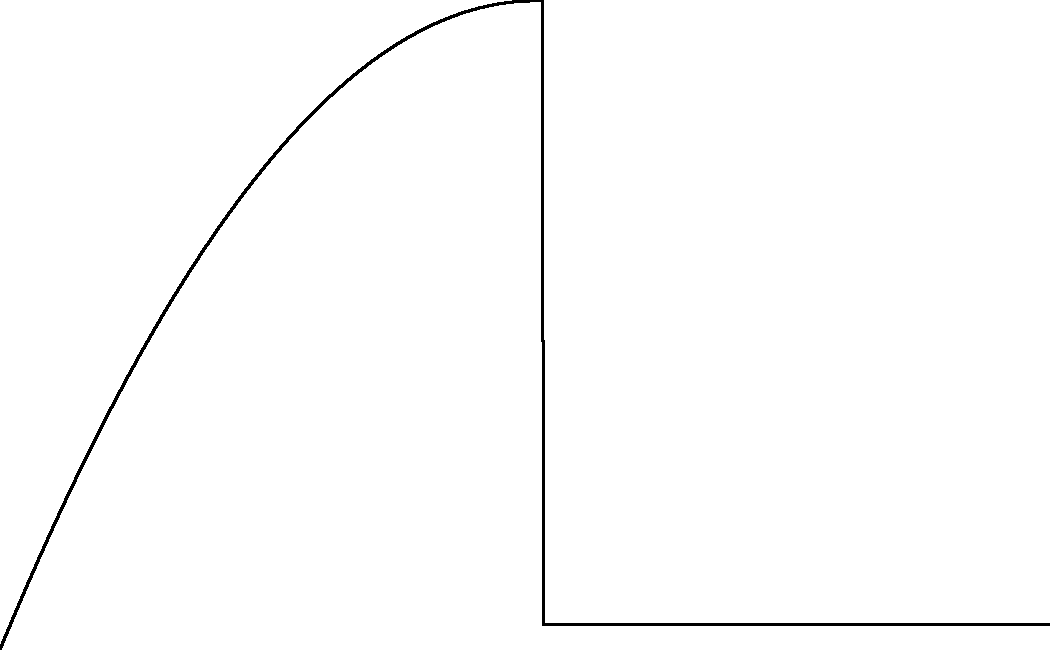}
	\end{minipage}%
	\begin{minipage}[t]{0.3\textwidth}
		\centering
		
		\includegraphics[width=0.9\textwidth]{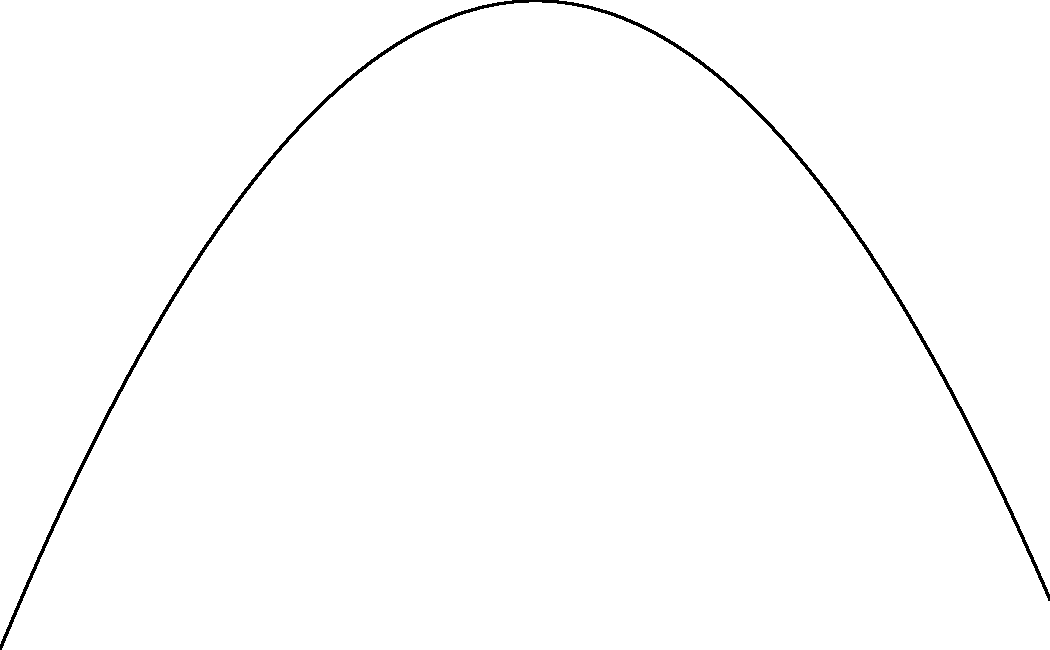}
	\end{minipage}
	\begin{minipage}[t]{0.3\textwidth}
		\centering
		
		\includegraphics[width=0.9\textwidth]{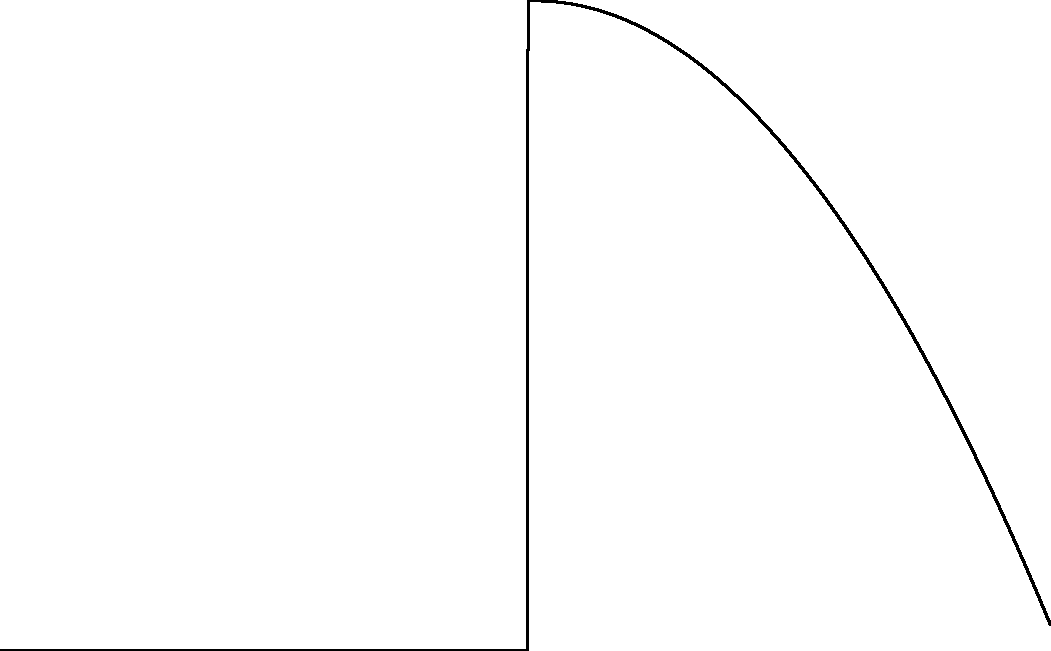}
	\end{minipage}
	\vskip0.5cm
	
	\begin{minipage}[t]{0.3\textwidth}
		\centering
		
		\includegraphics[width=0.9\textwidth]{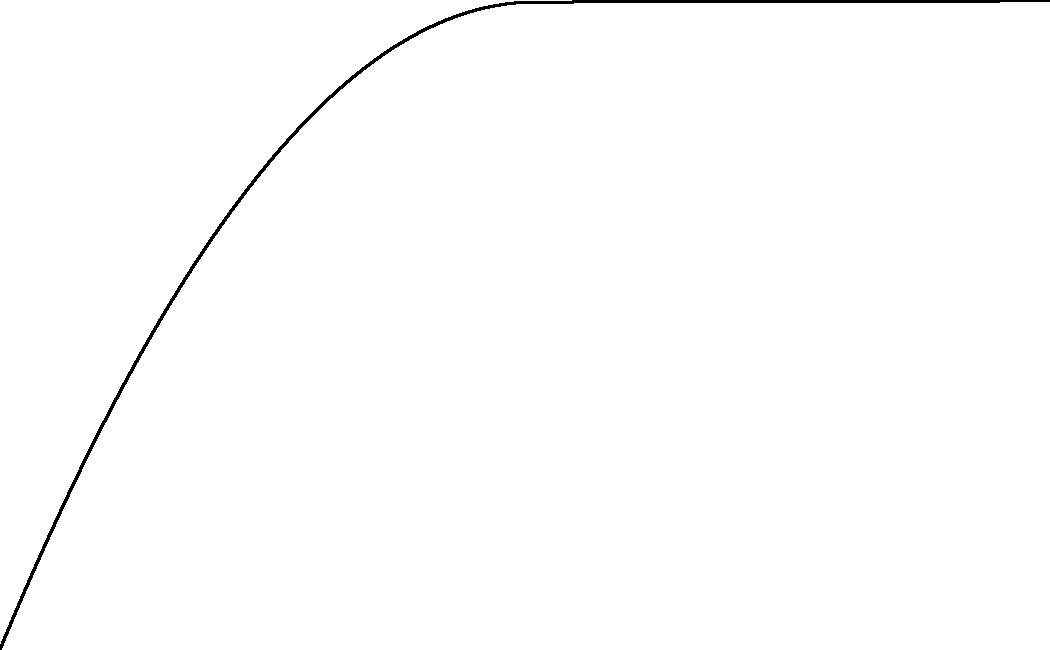}
	\end{minipage}%
	\begin{minipage}[t]{0.3\textwidth}
		\centering
		
		\includegraphics[width=0.9\textwidth]{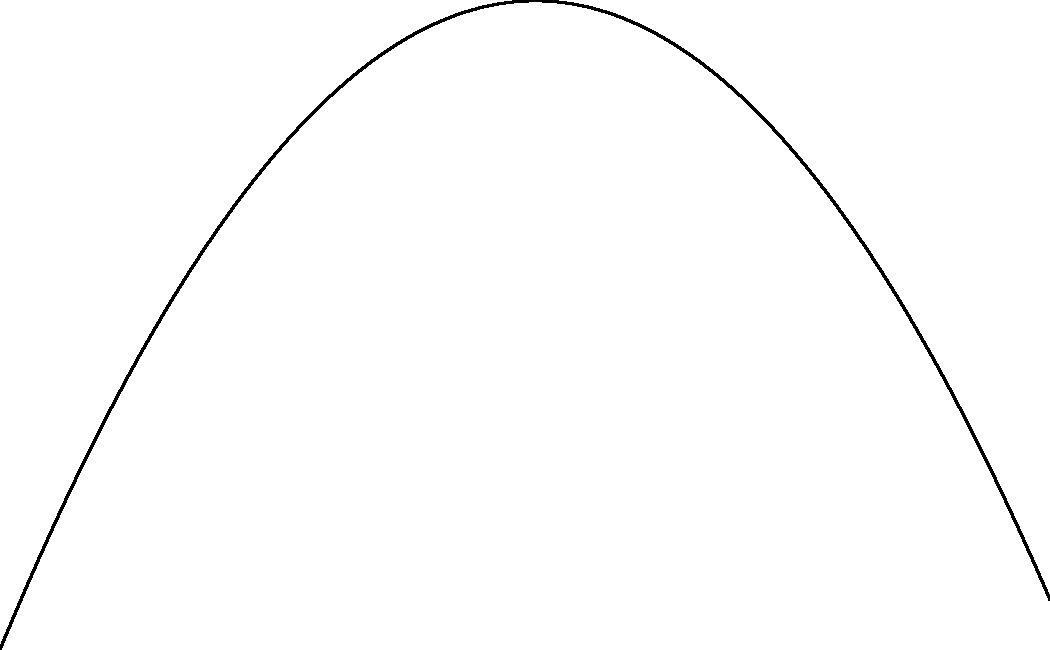}
	\end{minipage}
	\begin{minipage}[t]{0.3\textwidth}
		\centering
		
		\includegraphics[width=0.9\textwidth]{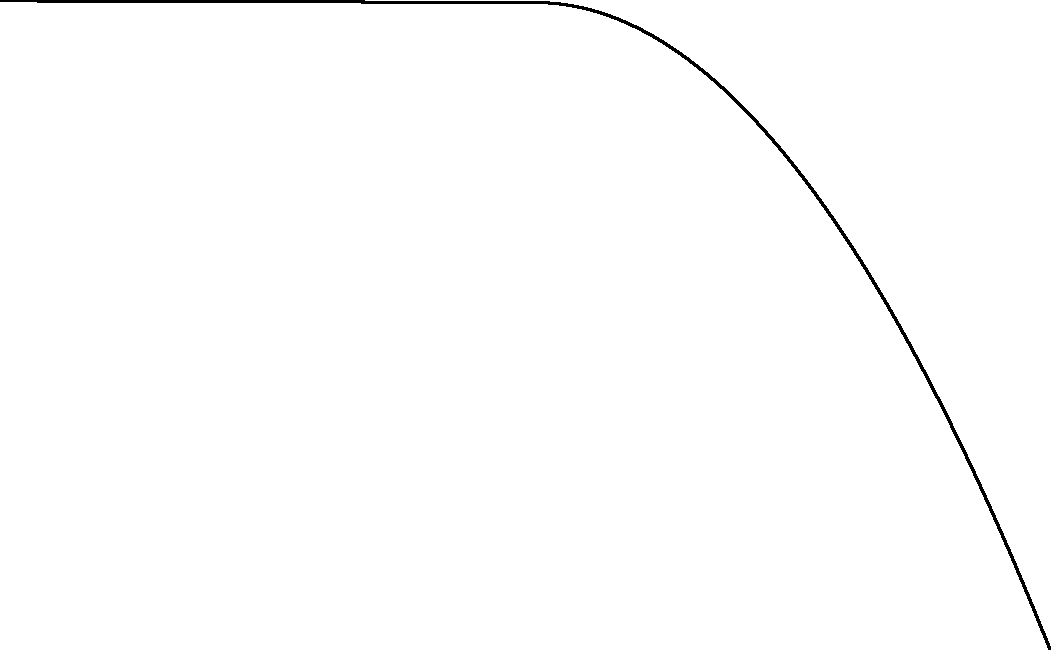}
	\end{minipage}
	\vskip0.5cm
	\caption{Numerical solutions to the smooth quasi-1D test problem (uniform mesh, $h=\frac{1}{1024}$). First row: sharp, $\delta=0$. Second row: sharp, $\delta=6h$. Third row: diffuse, $\delta=6h$. Fourth row: diffuse, $\delta=d_\Omega$.}
	\label{fig:s1d}
\end{figure}

Since the above exact solution is smooth and independent of $y$, we call this numerical experiment the smooth quasi-1D test. We performed grid convergence studies for both sharp and diffuse interface versions of our unfitted Nitsche method with PG stabilization. In particular, we studied the sensitivity of numerical results to the subdomains overlap $\delta$. The sharp interface simulations were run with $\delta = 0$ and $\delta = 6h$. The diffuse interface method was tested for $\delta = 6h$ and $\delta=d_\Omega$, where $d_\Omega = \mathrm{diam}(\bar{\Omega})$. For comparison purposes, we ran the same simulation with the Hansbo \& Hansbo ($\mathrm{H}^2$) method \cite{hansbo2}, which is essentially the FE formulation~\eqref{nitschewf}. The $L^2$ errors and experimental orders of convergence (EOC) on uniform grids are listed in Table~\ref{tab:1dsmooth}. In addition, we tested the sharp interface version with $\delta=6h$ on successively refined quasi-uniform grids, the coarsest of which is shown in Fig.~\ref{fig:coarse}. The $L^2$ error behavior is summarized in Table~\ref{tab:1dquasiu}. In all cases, we observe approximately second order convergence and (almost) no dependence on the value of $\delta$. The numerical solutions obtained on the uniform grid with $h=\frac{1}{1024}$ are shown in Fig.~\ref{fig:s1d}. As expected, a smooth transition from $u_1$ to $u_2$ can be observed across the interface.
\begin{table}[h!]
	\centering
	\begin{tabular}{ccccccccccc}
		\hline
		$h^{-1}$ & $\mathrm{H}^2$ & EOC & sharp & EOC  & sharp & EOC & diffuse & EOC  & diffuse & EOC\\
		& & &  $\delta=0$ & & $\delta=6h$ & & $\delta=0$ & & $\delta=d_\Omega$ &
		\\
		\hline 
		128  & 1.03e-05 &      & 4.02e-05 &      & 4.02e-05 &      & 4.02e-05 &      & 4.02e-05 &      \\ 
		256  & 2.59e-06 & 1.99 & 1.01e-05 & 1.99 & 1.01e-05 & 1.99 & 1.01e-05 & 1.99 & 1.01e-05 & 1.99 \\
		512  & 6.48e-07 & 2.00 & 2.54e-06 & 1.99 & 2.54e-06 & 1.99 & 2.54e-06 & 1.99 & 2.54e-06 & 1.99 \\
		1024 & 1.62e-07 & 2.00 & 6.35e-07 & 2.00 & 6.35e-07 & 2.00 & 6.35e-07 & 2.00 & 6.35e-07 & 2.00 \\ 
		2048 & 4.05e-08 & 2.00 & 1.59e-07 & 2.00 & 1.59e-07 & 2.00 & 1.59e-07 & 2.00 & 1.59e-07 & 2.00 \\
		4096 & 9.55e-09 & 2.08 & 3.96e-08 & 2.01 & 3.99e-08 & 1.99 & 3.98e-08 & 2.00 & 3.99e-08 & 1.99 \\
		\hline
	\end{tabular}
	\caption{Smooth quasi-1D test, $L^2$ convergence history on uniform meshes.}
	\label{tab:1dsmooth}
\end{table}

\subsection{Poisson problem with a straight interface and a non-smooth solution}
\label{sec:test1dkink}

In our second example, the domain and interface are the same as in the quasi-1D test of Section~\ref{sec:test1dsmooth}. Following Hansbo and Hansbo \cite{hansbo2}, we set $\mu_1=0.5$, $\mu_2=3$ and $f_1=f_2=1$. The Dirichlet boundary conditions at points $(x,y)\in\partial\Omega$ with $x=0$ or $x=1$ are chosen to match the values of
\begin{align*}
	u_1(x,y)= \frac{9}{14}(x-0.01)-(x-0.01)^2,\qquad u_2(x,y) = \frac{5}{84} + \frac{9}{84}(x-0.01) - \frac{1}{6}(x-0.01)^2.
\end{align*}
As before, we impose homogenous Neumann boundary conditions at $(x,y)\in\partial\Omega$ with $y=0$ or $y=1$.

\begin{figure}[h!]
	\centering
	\begin{minipage}[t]{0.3\textwidth}
		\centering (a)  extended $u_1$\vskip0.25cm
		
		\includegraphics[width=0.9\textwidth]{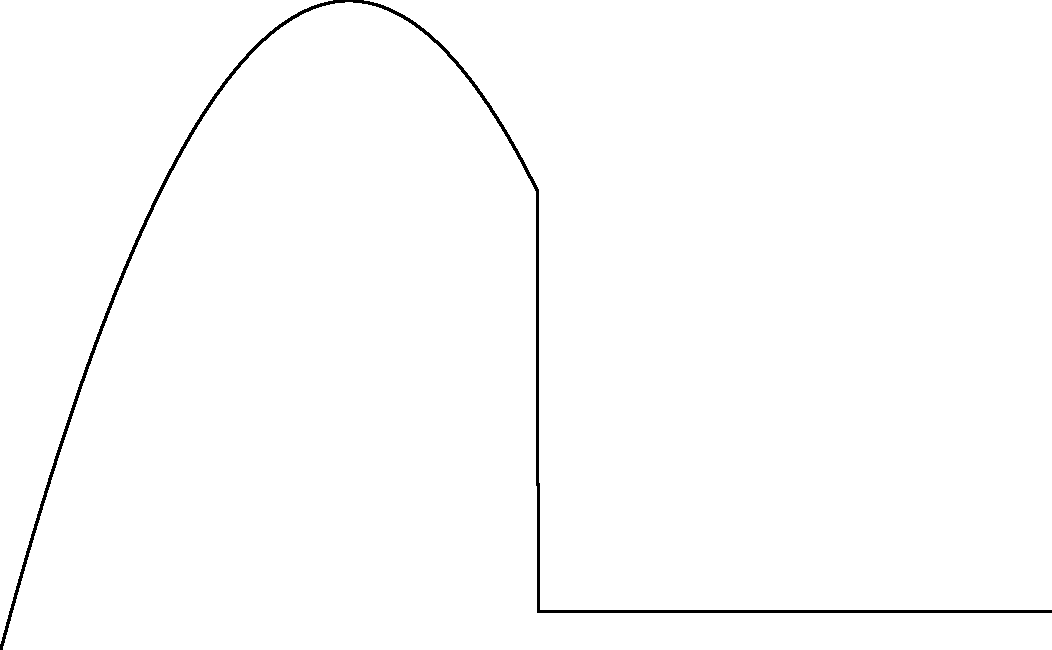}
	\end{minipage}%
	\begin{minipage}[t]{0.3\textwidth}
		\centering (b)  $u$\vskip0.25cm
		
		\includegraphics[width=0.9\textwidth]{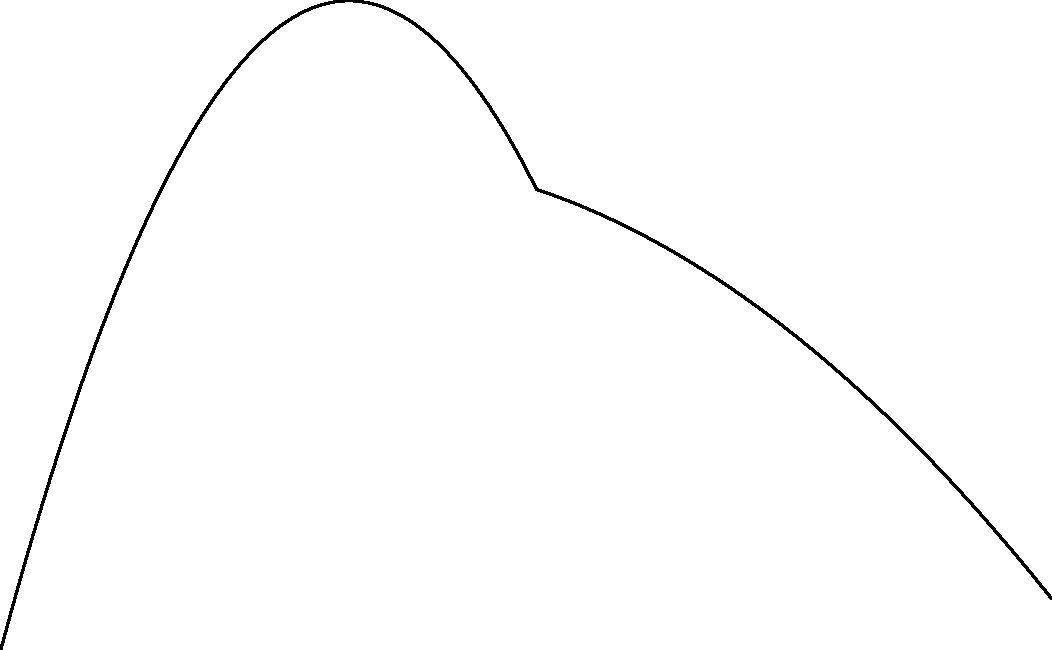}
	\end{minipage}
	\begin{minipage}[t]{0.3\textwidth}
		\centering (c)  extended $u_2$\vskip0.25cm
		
		\includegraphics[width=0.9\textwidth]{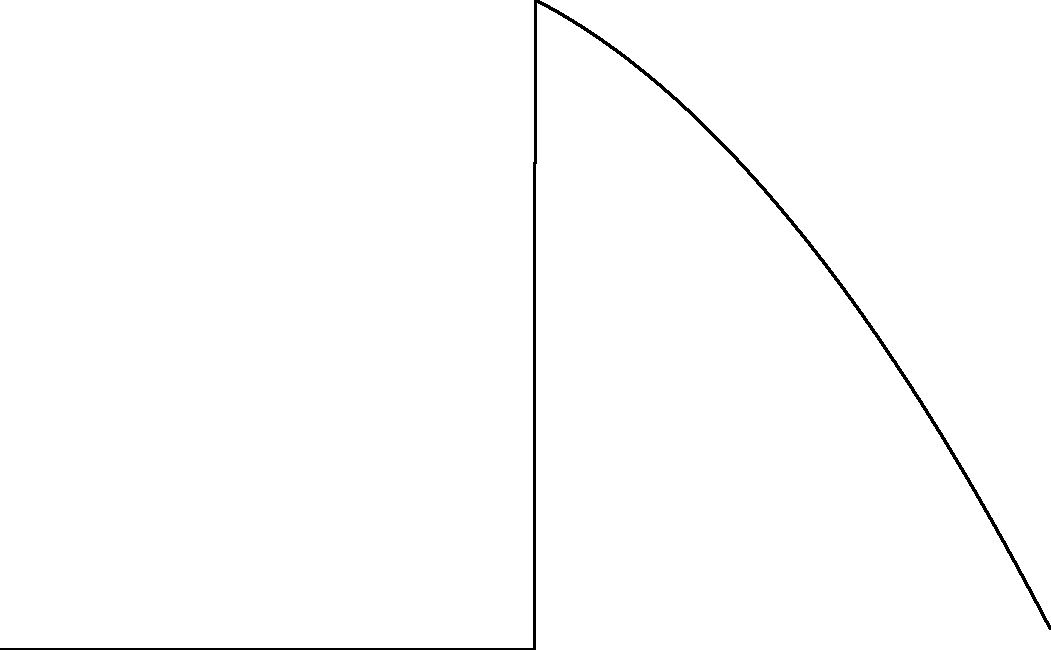}
	\end{minipage}
	\vskip0.5cm
	
	\begin{minipage}[t]{0.3\textwidth}
		\centering
		
		\includegraphics[width=0.9\textwidth]{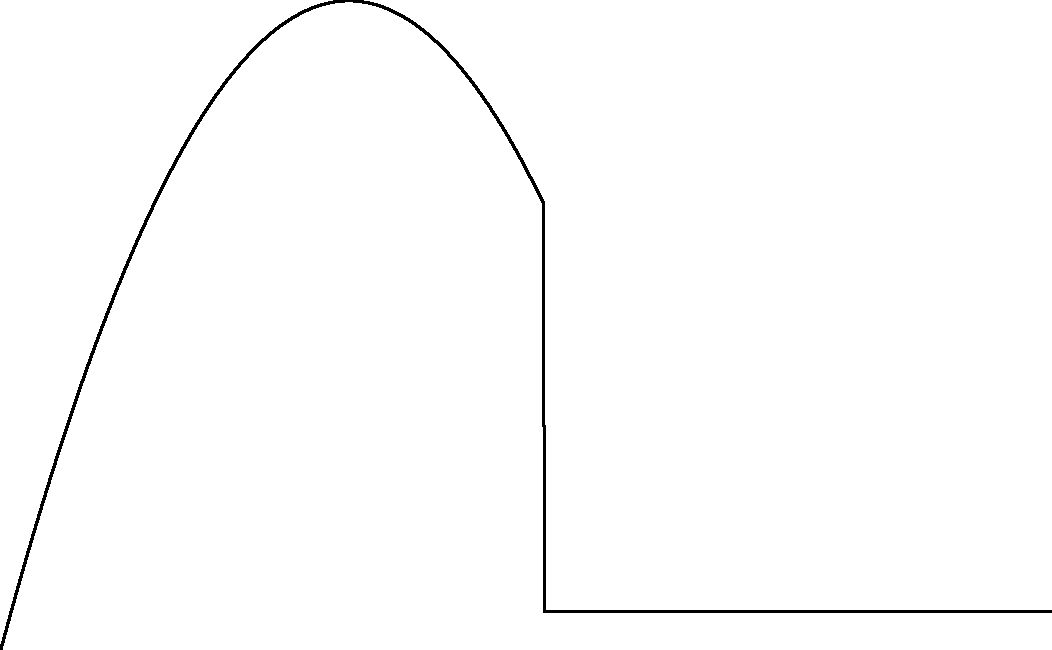}
	\end{minipage}%
	\begin{minipage}[t]{0.3\textwidth}
		\centering
		
		\includegraphics[width=0.9\textwidth]{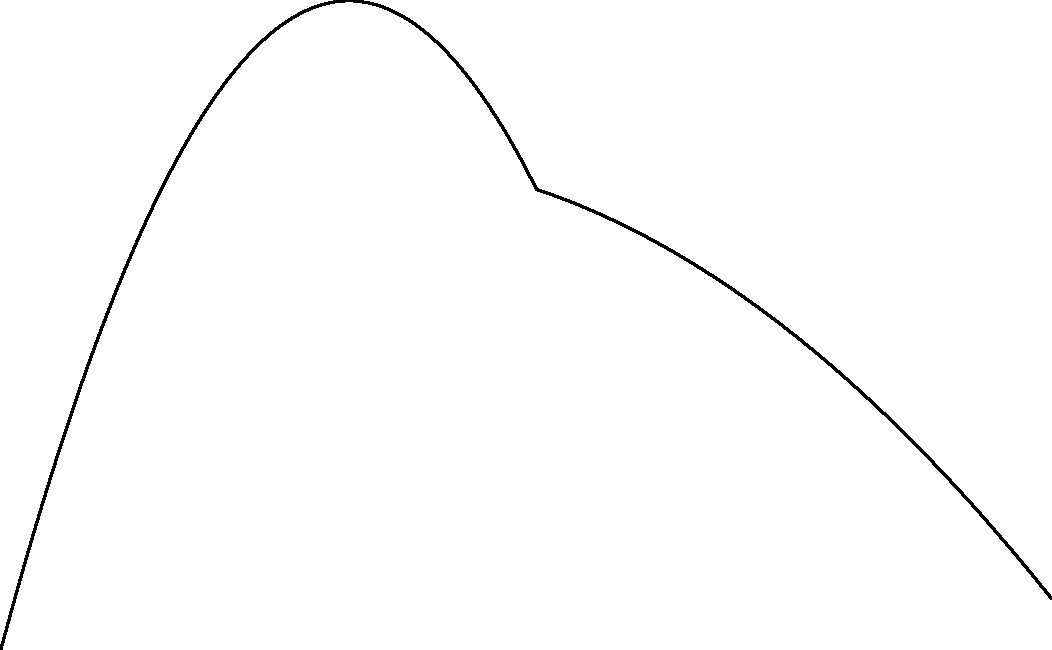}
	\end{minipage}
	\begin{minipage}[t]{0.3\textwidth}
		\centering
		
		\includegraphics[width=0.9\textwidth]{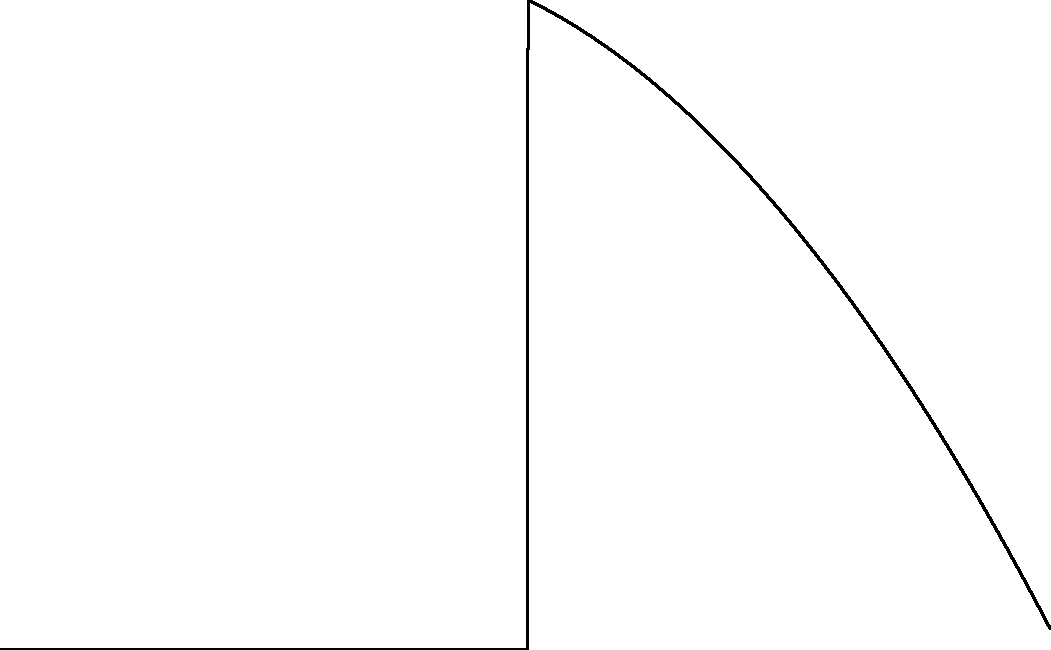}
	\end{minipage}
	\vskip0.5cm
	
	\begin{minipage}[t]{0.3\textwidth}
		\centering
		
		\includegraphics[width=0.9\textwidth]{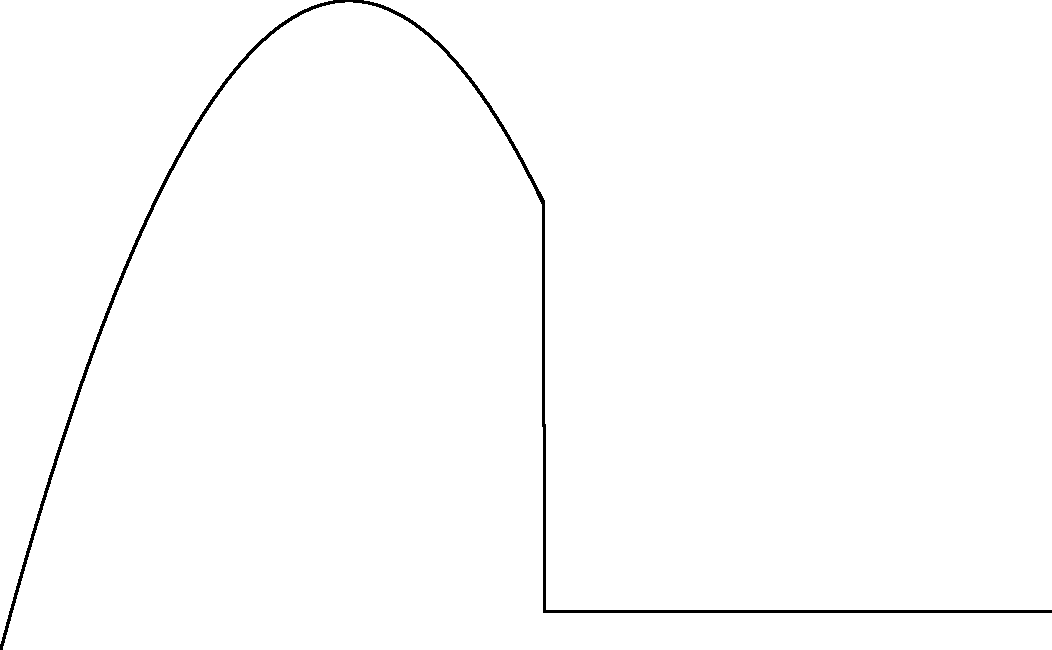}
	\end{minipage}%
	\begin{minipage}[t]{0.3\textwidth}
		\centering
		
		\includegraphics[width=0.9\textwidth]{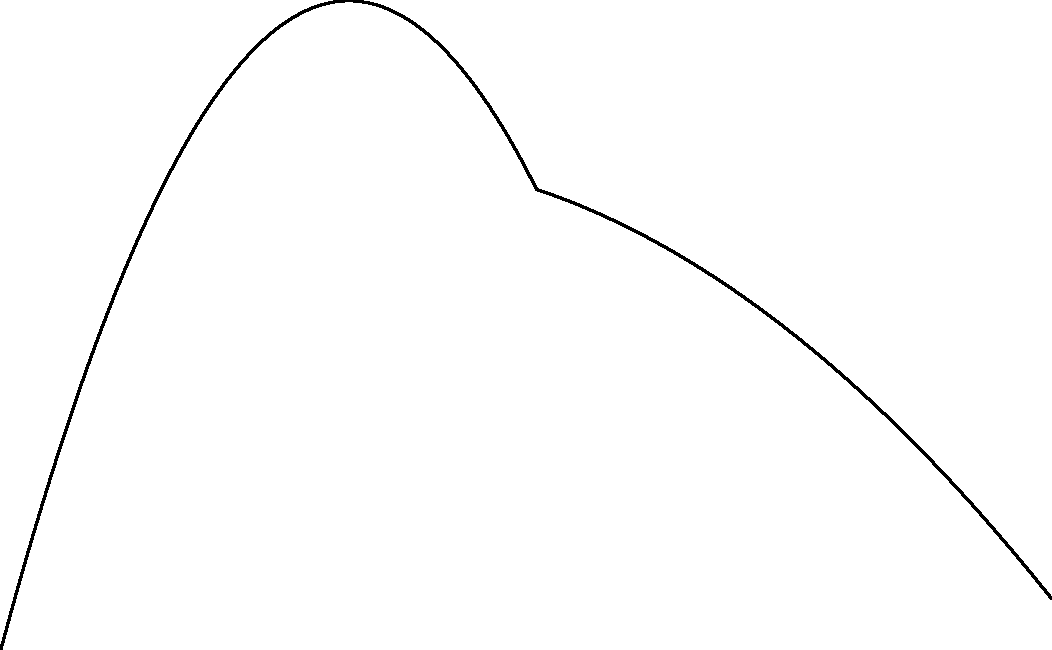}
	\end{minipage}
	\begin{minipage}[t]{0.3\textwidth}
		\centering
		
		\includegraphics[width=0.9\textwidth]{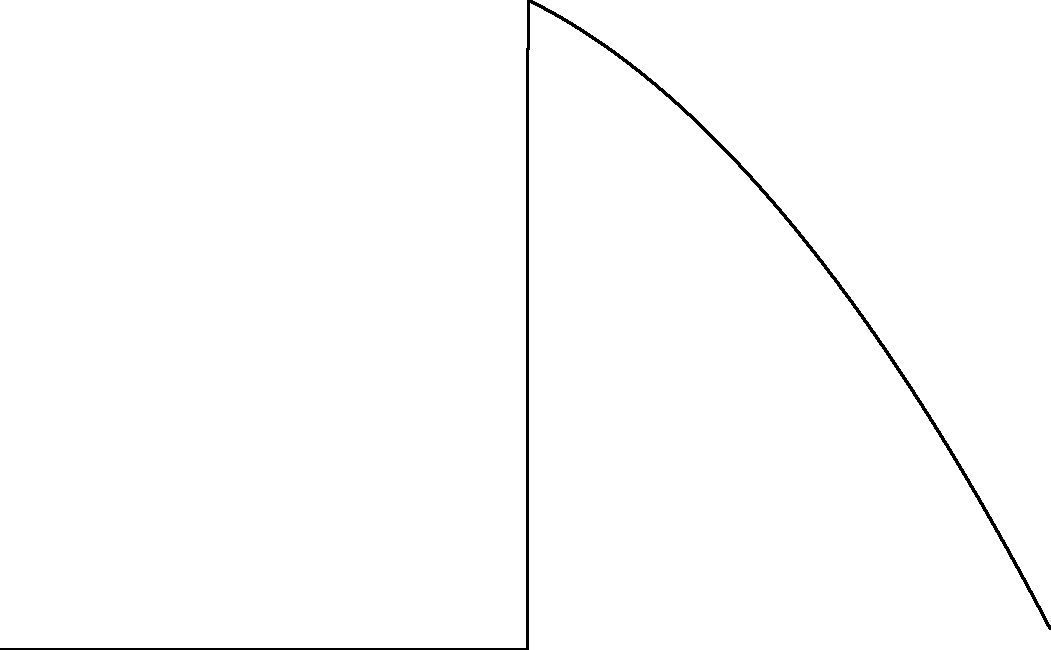}
	\end{minipage}
	\vskip0.5cm
	
	\begin{minipage}[t]{0.3\textwidth}
		\centering
		
		\includegraphics[width=0.9\textwidth]{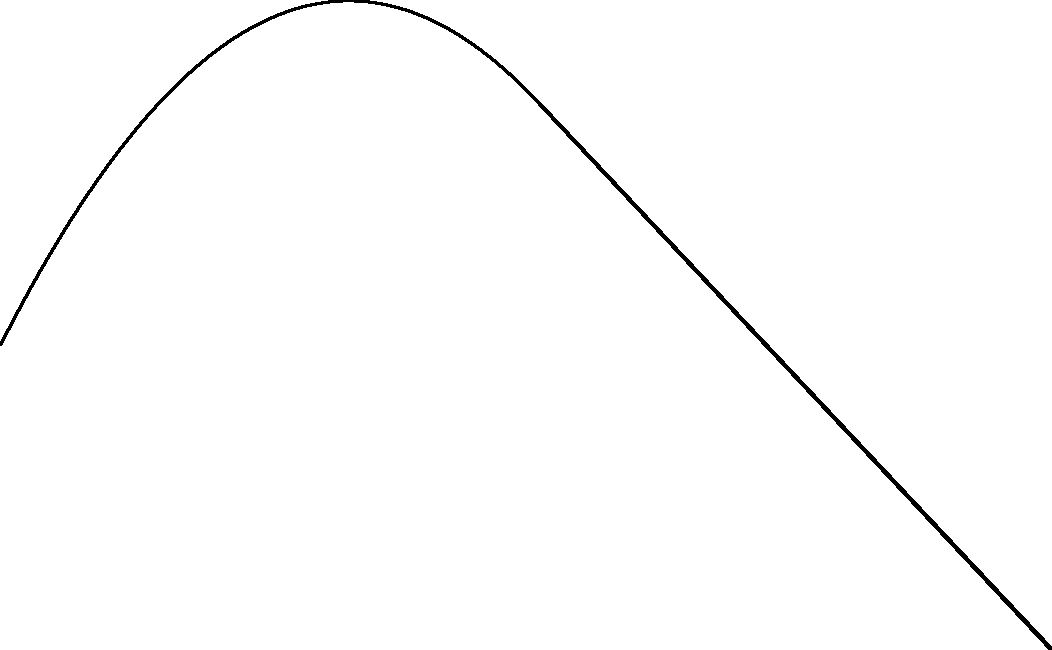}
	\end{minipage}%
	\begin{minipage}[t]{0.3\textwidth}
		\centering
		
		\includegraphics[width=0.9\textwidth]{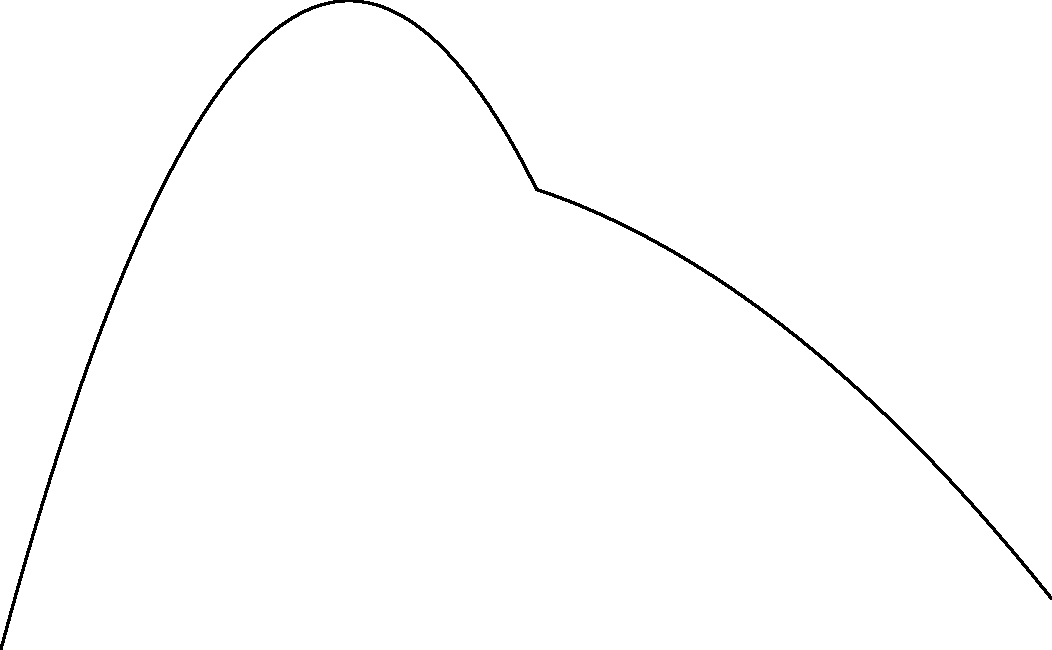}
	\end{minipage}
	\begin{minipage}[t]{0.3\textwidth}
		\centering
		
		\includegraphics[width=0.9\textwidth]{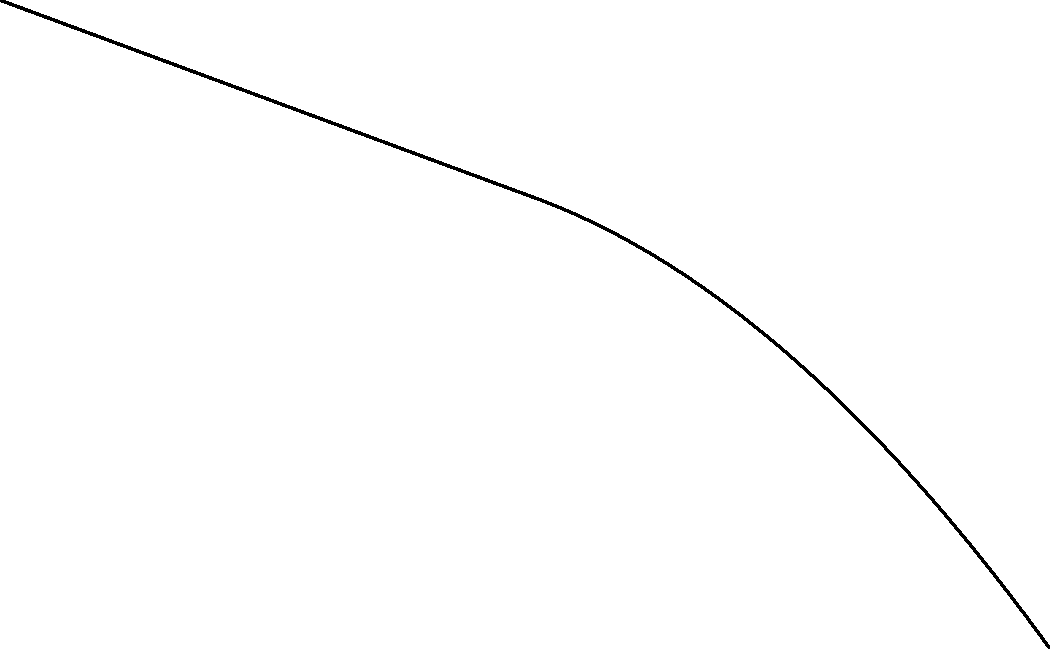}
	\end{minipage}
	\vskip0.5cm
	\caption{Numerical solutions to the non-smooth quasi-1D test problem (uniform mesh, $h=\frac{1}{256}$). First row: sharp, $\delta=0$. Second row: sharp, $\delta=6h$. Third row: diffuse, $\delta=6h$. Fourth row: diffuse, $\delta=d_\Omega$.}
	\label{fig:k1d}
\end{figure}

The analytical solution $u(x,y)$ to this test problem coincides with $u_1(x,y)$ for $x\le 0.51$ and with $u_2(x,y)$ for $x\ge 0.51$. Since $u(x,y)$ has a kink at $x=0.51$, we call this experiment the non-smooth quasi-1D test. We ran grid convergence studies for the same configurations as in the previous example. The $L^2$ errors and EOCs for uniform meshes are listed in Table \ref{tab:1dkink}. Additionally, we present the $L^2$ convergence history for quasi-uniform meshes in Table ~\ref{tab:1dquasiu}. Again, we observe second order convergence and no dependence on $\delta$ for all cases. The numerical solutions obtained with $h=\frac{1}{1024}$ are shown in Fig.~\ref{fig:k1d}. The solutions for $u_1$ and $u_2$ match well on the interface and no oscillations occur.

\begin{rmrk}
	The reported convergence behavior can be further improved by using problem-dependent values of the Nitsche parameter $\alpha$. Numerical tests indicate that more stable EOCs can often be obtained by increasing or decreasing the value of $\alpha$. 
\end{rmrk}

\begin{table}[h!]
	\centering
	\begin{tabular}{ccccc}
		\hline
		$h^{-1}$ & smooth & EOC & non-smooth & EOC\\
		\hline
		10  & 4.23e-03 &      & 3.09e-03 &      \\
		20  & 1.15e-03 & 1.88 & 8.36e-04 & 1.89 \\
		40  & 2.98e-04 & 1.94 & 2.16e-04 & 1.95 \\
		80  & 7.60e-05 & 1.98 & 5.50e-05 & 1.97 \\
		160 & 1.92e-05 & 1.98 & 1.39e-05 & 1.98 \\
		320 & 4.82e-06 & 1.99 & 3.48e-06 & 2.00 \\
		640 & 1.21e-06 & 1.99 & 8.72e-07 & 2.00 \\
		\hline
	\end{tabular}
	\caption{Quasi-1D test problems, $L^2$ convergence of the sharp interface method on quasi-uniform meshes.}
	\label{tab:1dquasiu}
\end{table}

\begin{table}[h!]
	\centering
	\begin{tabular}{ccccccccccc}
		\hline
		$h^{-1}$ & $\mathrm{H}^2$ & EOC & sharp & EOC  & sharp & EOC & diffuse & EOC  & diffuse & EOC\\
		& & &  $\delta=0$ & & $\delta=6h$ & & $\delta=6h$ & & $\delta=d_\Omega$ &
		\\
		\hline 
		128  & 7.47e-06 &      & 2.91e-05 &      & 2.91e-05 &      & 2.67e-05 &      & 2.67e-05 &      \\ 
		256  & 1.86e-06 & 2.01 & 7.31e-06 & 1.99 & 7.31e-06 & 1.99 & 8.29e-06 & 1.69 & 8.29e-06 & 1.69 \\
		512  & 4.64e-07 & 2.00 & 1.83e-06 & 2.00 & 1.83e-06 & 2.00 & 1.38e-06 & 2.59 & 1.38e-06 & 2.59 \\
		1024 & 1.15e-07 & 2.01 & 4.57e-07 & 2.00 & 4.57e-07 & 2.00 & 4.16e-07 & 1.73 & 4.16e-07 & 1.73 \\ 
		2048 & 1.93e-08 & 2.57 & 1.03e-07 & 2.15 & 1.03e-07 & 2.15 & 9.79e-08 & 2.09 & 9.79e-07 & 2.09 \\
		4096 & 4.64e-09 & 2.06 & 2.29e-08 & 2.17 & 2.29e-08 & 2.17 & 2.75e-08 & 1.83 & 2.75e-08 & 1.83 \\
		\hline
	\end{tabular}
	\caption{Non-smooth quasi-1D test, $L^2$ convergence history on uniform meshes.}
	\label{tab:1dkink}
\end{table}

\subsection{Poisson problem with a circular interface}
Let us now apply our PG method to a truly two-dimensional test problem with a circular interface. The domain is given by $\Omega=(-1,1)^2$. The interface $\Gamma$ is the zero level set of the signed distance function $\phi(x,y)=0.75-\sqrt{x^2+y^2}$. We set $\mu_1=1$, $\mu_2=10^3$ and $f_1=f_2=4$ in this test. The Dirichlet boundary condition for points $(x,y)\in\partial\Omega$ is chosen in such a way that the restrictions
\begin{align*}
	u_1(x,y) =x^2-y^2,
	\qquad u_2(x,y) = \frac{x^2-y^2}{1000}-\frac{0.5625}{1000}+0.5625
\end{align*}
define the analytical solution (cf. \cite{hansbo2}). 

We studied grid convergence of the sharp and diffuse interface version of our unfitted Nitsche method with PG stabilization. The parameter settings were chosen as in the previous examples. For comparison purposes, the same simulations were run with the $\mathrm{H}^2$ method. The $L^2$ errors and EOCs are listed in Table \ref{tab:2d}. We observed second order convergence in all cases.

\begin{figure}[h!]
	\centering
	\begin{minipage}[t]{0.33\textwidth}
		\centering (a)  extended $u_1$\vskip0.25cm
		
		\includegraphics[width=0.9\textwidth]{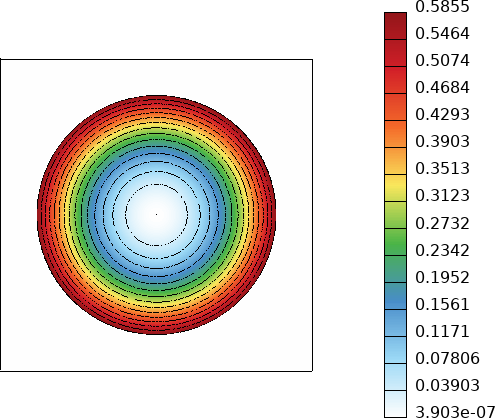}
	\end{minipage}%
	\begin{minipage}[t]{0.33\textwidth}
		\centering (b)  $u$\vskip0.25cm
		
		\includegraphics[width=0.9\textwidth]{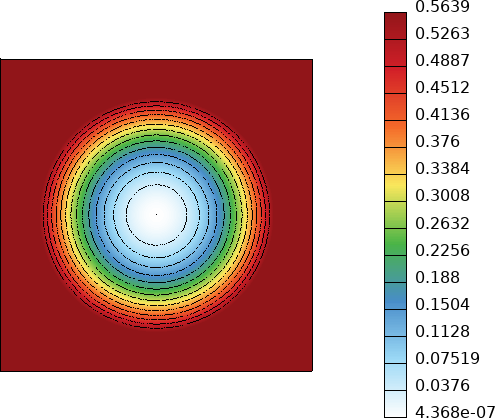}
	\end{minipage}
	\begin{minipage}[t]{0.33\textwidth}
		\centering (c)  extended $u_2$\vskip0.25cm
		
		\includegraphics[width=0.9\textwidth]{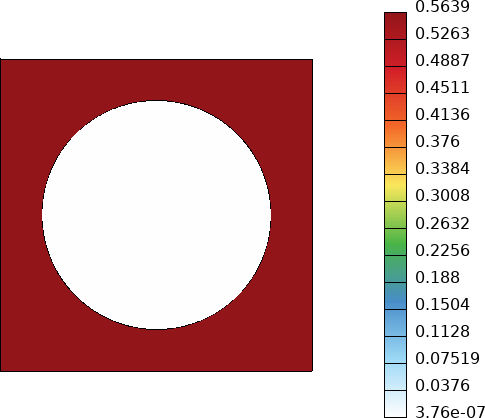}
	\end{minipage}
	\vskip0.5cm
	
	\begin{minipage}[t]{0.33\textwidth}
		\includegraphics[width=0.9\textwidth]{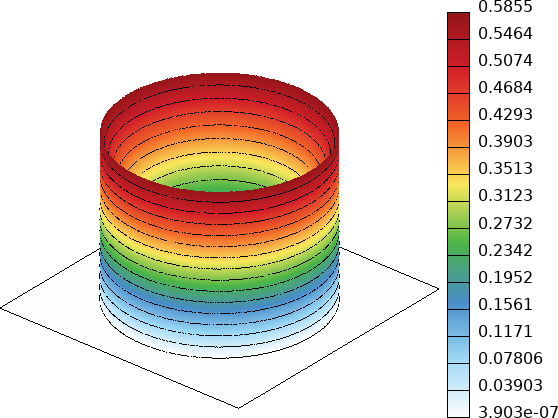}
	\end{minipage}%
	\begin{minipage}[t]{0.33\textwidth}
		\includegraphics[width=0.9\textwidth]{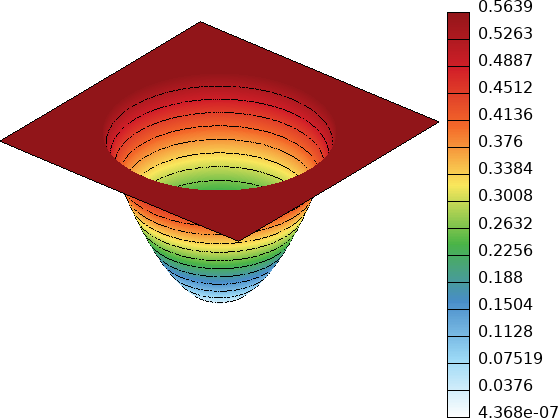}
	\end{minipage}
	\begin{minipage}[t]{0.33\textwidth}
		\includegraphics[width=0.9\textwidth]{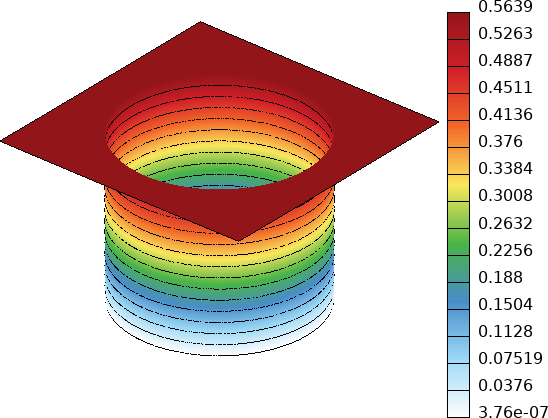}
	\end{minipage}
	\vskip0.25cm 
	\caption{Numerical solutions to the circular interface test problem (uniform mesh, $h=\frac{1}{512}$, diffuse interface method, $\delta=6h$). First row: top views. Second row: elevation is proportional to the value of the plotted function.}
	\label{fig:k2a}
\end{figure}

\begin{figure}[h!]
	\centering
	\begin{minipage}[t]{0.33\textwidth}
		\centering (a)  extended $u_1$\vskip0.5cm
		
		\includegraphics[width=0.9\textwidth]{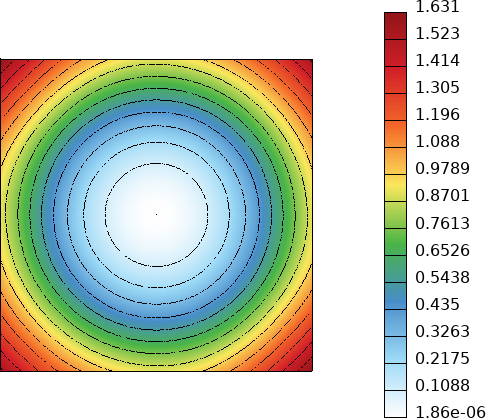}
	\end{minipage}%
	\begin{minipage}[t]{0.33\textwidth}
		\centering (b)  $u$\vskip0.5cm
		
		\includegraphics[width=0.9\textwidth]{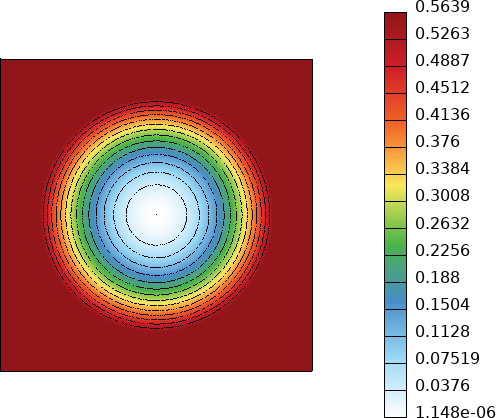}
	\end{minipage}
	\begin{minipage}[t]{0.33\textwidth}
		\centering (c)  extended $u_2$\vskip0.5cm
		
		\includegraphics[width=0.9\textwidth]{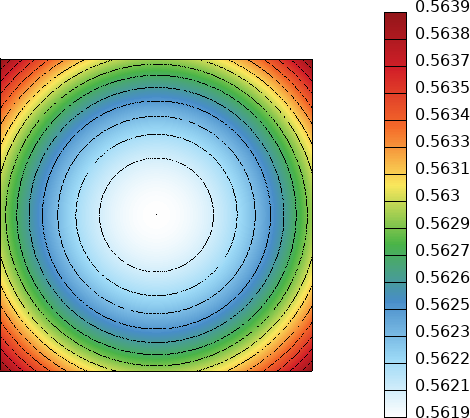}
	\end{minipage}
	\vskip0.5cm
	
	\begin{minipage}[t]{0.33\textwidth}
		\includegraphics[width=0.9\textwidth]{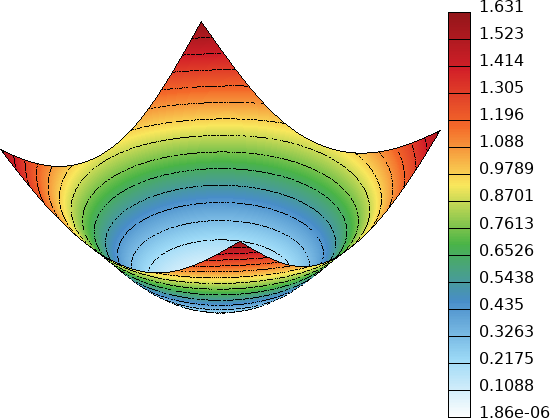}
	\end{minipage}%
	\begin{minipage}[t]{0.33\textwidth}
		\includegraphics[width=0.9\textwidth]{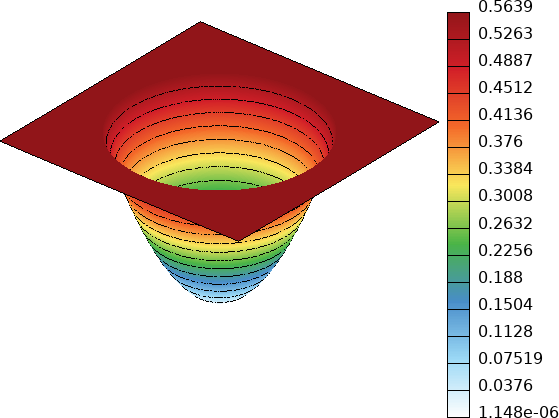}
	\end{minipage}
	\begin{minipage}[t]{0.33\textwidth}
		\includegraphics[width=0.9\textwidth]{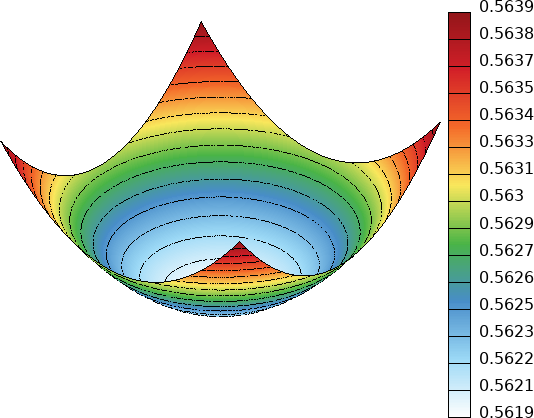}
	\end{minipage}
	\vskip0.25cm 
	\caption{Numerical solutions to the circular interface test problem (uniform mesh, $h=\frac{1}{512}$, diffuse interface method, $\delta = d_\Omega$). First row: top views. Second row: elevation is proportional to the value of the plotted function.}
	\label{fig:k2b}
\end{figure}

The numerical solutions corresponding to $h=\frac{1}{512}$ are shown in Figs.~\ref{fig:k2a} and \ref{fig:k2b} for diffuse interface simulations with $\delta=6h$ and $\delta = d_\Omega$, respectively. It can be seen that the unfitted finite element approximations to $u_1$ and $u_2$ match well on the circular interface. Moreover, no violations of the discrete maximum principle are observed in this experiment. 

\begin{table}[h!]
	\centering
	\rev{
	\begin{tabular}{ccccccccccc}
		\hline
		$h^{-1}$ & $\mathrm{H}^2$ & EOC & sharp & EOC  & sharp & EOC & diffuse & EOC  & diffuse & EOC\\
		& & &  $\delta=0$ & & $\delta=6h$ & & $\delta=6h$ & & $\delta=d_\Omega$ & \\
		\hline
		8    & 2.79e-02 &      & 6.33e-02 &      & 7.04e-02 &      & 7.03e-02 &      & 7.03e-02 &      \\
		16   & 6.89e-03 & 2.02 & 1.69e-02 & 1.91 & 1.97e-02 & 1.84 & 2.01e-02 & 1.81 & 2.01e-02 & 1.81 \\
		32   & 1.74e-03 & 1.99 & 4.32e-03 & 1.97 & 5.16e-03 & 1.93 & 5.14e-03 & 1.97 & 5.16e-03 & 1.96 \\
		64   & 4.33e-04 & 2.01 & 1.09e-03 & 1.99 & 1.25e-03 & 2.05 & 1.24e-03 & 2.05 & 1.32e-03 & 1.97 \\
		128  & 1.08e-04 & 2.00 & 2.74e-04 & 1.99 & 2.96e-04 & 2.08 & 2.90e-04 & 2.10 & 3.30e-04 & 2.00 \\
		256  & 2.69e-05 & 2.01 & 6.87e-05 & 2.00 & 7.16e-05 & 2.05 & 6.86e-05 & 2.08 & 8.15e-05 & 2.02 \\
		512  & 6.72e-06 & 2.00 & 1.72e-05 & 2.00 & 1.76e-05 & 2.02 & 1.61e-05 & 2.09 & 1.97e-05 & 2.05 \\
		1024 & 1.69e-06 & 1.99 & 4.31e-06 & 2.00 & 4.35e-06 & 2.02 & 3.60e-06 & 2.16 & 4.56e-06 & 2.11 \\
		\hline
	\end{tabular}}
	\caption{Circular interface test, $L^2$ convergence history on uniform meshes.}
	\label{tab:2d}
\end{table}

\rev{\subsection{High-order finite elements and condition number of the system matrix}
\label{sec:HO}

Let us now demonstrate numerically  that optimal-order accuracy is attained with PG stabilization for higher order finite elements. In our numerical example, the computational domain is given by $\Omega = (0,1)^2$. As in Sections~\ref{sec:test1dsmooth} and~\ref{sec:test1dkink}, the interface is located at $x = 0.51$. The diffusion coefficients are set to $\mu_1 = 0.5$ and $\mu_2 = 2$. The right-hand side is constructed from an exact solution satisfying
\begin{equation*}
u_1(x,y) = \frac{7}{12}(x - 0.51) - \frac{2}{3}(x - 0.51)^4,\qquad
u_2(x,y) = \frac{25}{176} + \frac{47}{176}(x - 0.51) - \frac{9}{22}(x - 0.51)^4
\end{equation*}
in the subdomains $\Omega_1$ and $\Omega_2$, respectively. The Dirichlet boundary conditions imposed at $(x,y)\in\partial\Omega$ with $x = 0$ or $x = 1$ are chosen to match the exact solution. Homogeneous Neumann conditions are imposed at $(x,y)\in\partial\Omega$ with $y = 0$ or $y = 1$.

We investigate the grid convergence behavior of our sharp interface method with $\delta = 0$ and again compare it against the $\mathrm{H}^2$ method. The $L^2$ errors and experimental orders of convergence (EOCs) for polynomial degrees $p = 1,2,3$ are reported in Table~\ref{tab:hohansbo} for the $\mathrm{H}^2$ method and in Table~\ref{tab:hostab} for the sharp interface version of our unfitted Nitsche method with PG stabilization. In all cases, the observed convergence rate is $p+1$.

\begin{table}[h!]
	\centering
	\rev{
	\begin{tabular}{ccccccc}
		\hline
		$h^{-1}$ & $\mathrm{H}^2$ & EOC & $\mathrm{H}^2$ & EOC  & $\mathrm{H}^2$ & EOC\\
		& $p=1$ & & $p=2$ & & $p=3$ & \\
		\hline 
		32   & 1.67e-04 &      & 9.79e-07 &      & 3.60e-09 &      \\ 
		64   & 4.19e-05 & 1.99 & 1.24e-07 & 2.98 & 2.26e-10 & 3.99 \\
		128  & 1.05e-05 & 2.00 & 1.56e-08 & 2.99 & 1.42e-11 & 3.99 \\
		256  & 2.63e-06 & 2.00 & 1.95e-09 & 3.00 &          &      \\ 
		512  & 6.57e-07 & 2.00 &          &      &          &      \\
		\hline
	\end{tabular}
	\caption{Quasi-1D test with quartic solution, $L^2$ convergence history of the $\mathrm{H}^2$ method on uniform meshes for different polynomial approximations.}
	\label{tab:hohansbo}}
\end{table}

\begin{table}[h!]
	\centering
	\rev{
	\begin{tabular}{ccccccc}
		\hline
		$h^{-1}$ & sharp & EOC & sharp & EOC  & sharp & EOC\\
		& $p=1$ & & $p=2$ & & $p=3$ & \\
		\hline 
		32   & 6.51e-04 &      & 9.92e-07 &      & 6.05e-09 &      \\ 
		64   & 1.70e-04 & 1.94 & 1.24e-07 & 3.00 & 3.86e-10 & 3.97 \\
		128  & 4.35e-05 & 1.97 & 1.56e-08 & 2.99 & 2.45e-11 & 3.99 \\
		256  & 1.10e-05 & 1.98 & 1.95e-09 & 3.00 &          &      \\ 
		512  & 2.77e-06 & 1.99 &          &      &          &      \\
		\hline
	\end{tabular}
	\caption{Quasi-1D test with quartic solution, $L^2$ convergence history of our stabilized sharp interface method with $\delta = 0$ on uniform meshes for different polynomial approximations.}
	\label{tab:hostab}}
\end{table}

In addition to assessing the convergence rates, we investigate the dependence of the condition number on the interface position within the mesh for linear finite elements. Specifically, the interface is placed at $x = 0.5 + 10^{-j}$ for $j = 2, 3, \ldots, 9$ on a fixed mesh, and the condition number of the system matrix is evaluated for both the $\mathrm{H}^2$ method and the proposed unfitted sharp interface method with $\delta = 0$. The condition numbers, scaled by $h^2$, are reported in Table~\ref{tab:cond}. For the $\mathrm{H}^2$ method, the condition number increases significantly as the cut cells become smaller. For the stabilized sharp interface method, the condition number remains uniformly bounded, which demonstrates the algebraic stability of the proposed approach.

\begin{table}[h!]
	\centering
	\rev{
	\begin{tabular}{ccccccccc}
		\hline
		$\mathrm{dist}(\partial K,\Gamma)$ & 1e-09 & 1e-08 & 1e-07 & 1e-06  & 1e-05 & 1e-04 & 1e-03 & 1e-02 
		\\
		\hline 
		$\mathrm{H}^2$  & 9.28e14 & 5.48e13 & 5.56e11 & 5.56e09 & 5.51e07 & 5.05e05 & 2.81e03 & 16.37\\
		$\delta = 0$    & 22.32   & 22.32   & 22.32   & 22.32   & 22.31   & 22.23   & 21.49   & 15.99\\
		\hline
	\end{tabular}
	\caption{Scaled condition numbers  $\kappa h^2$ of the system matrix for different distances between the mesh cell boundary and interface.}
	\label{tab:cond}}
\end{table}
}

\subsection{Stationary convection-diffusion problem with steep gradients}
\label{sec:convection}

In our final example, we test a generalization of our diffuse interface method to the elliptic problem
\begin{align*}
	\mathbf{v}\cdot\nabla u-\nabla\cdot(\mu\nabla u)&=0\quad\text{in }\Omega_k,~{k=1,2},\label{eq:pde}\\
	u &= g_D\quad\text{on }\partial\Omega,\\
	[\![ u]\!]&=0\quad\text{on }\Gamma,\\
	[\![\mathbf v u-\mu\nabla u]\!]&=0\quad\text{on }\Gamma,
\end{align*}
where $\mathbf v\in\mathbb R^d$ is a constant velocity and $g_D$ is the
Dirichlet boundary data. In a corresponding straightforward
extension of our diffuse
interface formulation, we add the convective term
$$\mathbf v\cdot\left[
\int_{\Omega_1}\nabla u_{h,1}w_{h,1}\xdif \bx+
\int_{\Omega_2}\nabla u_{h,2}w_{h,2}\xdif \bx\right]
$$
on the left-hand side of \eqref{wdiab} and multiply the stabilization
term $s_h(u,w)$ by (cf. \cite{weno})  
$$
\beta_h=\max\left\{1,\frac{|\mathbf v| h}{2\mu}\right\}.
$$

A well-known test problem for stationary perturbed convection-diffusion
equations with constant coefficients is defined in Example 2 of the
review paper \cite{john3}. It uses the domain $\Omega=(0,1)^2$ and
the velocity $\mathbf v=(\cos(-\pi/3),\sin(-\pi/3))^\top$. The
Dirichlet boundary data is given by
$$
g_D(x,y)=\begin{cases}
	0 & \mbox{if}\ x=1\ \mbox{or}\ y\le 0.7,\\ 1 & \mbox{otherwise}.
\end{cases}
$$
Note that $g_D$ is discontinuous at $(x_0,y_0)=(0,0.7)$. Let $\Gamma
=\{(x,y)\in\bar\Omega\,: ax+by=c\}$
be the straight line that passes through $(x_0,y_0)$ and is parallel
to $\mathbf v$. Deviating from the setup in \cite{john3}, 
we set $\mu_1=10^{-3}$ in
$\Omega_1=\{(x,y)\in\Omega\,: ax+by<c\}$ and $\mu_2=10^{-8}$ in
$\Omega_2=\{(x,y)\in\Omega\,: ax+by>c\}$.

\rev{
\begin{remark}
We use a constant velocity in this experiment to match the standard benchmark problem from \cite[Example~2]{john3}. The effectiveness of the projected gradient stabilization for arbitrary velocity fields and even for high-order finite element discretizations of nonlinear hyperbolic systems has been demonstrated in \cite{entropyHO,lohmann2017,weno}.
\end{remark}
}

\begin{figure}[h!]
	\centering
	\begin{minipage}[t]{0.33\textwidth}
		\centering (a)  no stabilization, $\beta_h = 0$\vskip0.5cm
		
		\includegraphics[width=0.9\textwidth]{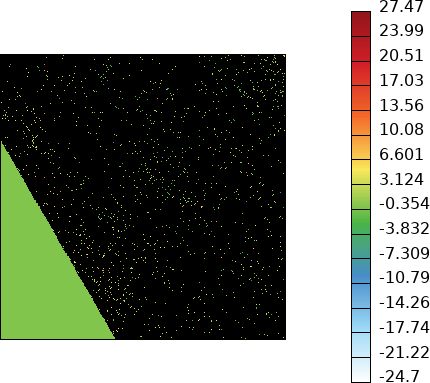}
	\end{minipage}%
	\begin{minipage}[t]{0.33\textwidth}
		\centering (b)  $\beta_h = 1$\vskip0.5cm
		
		\includegraphics[width=0.9\textwidth]{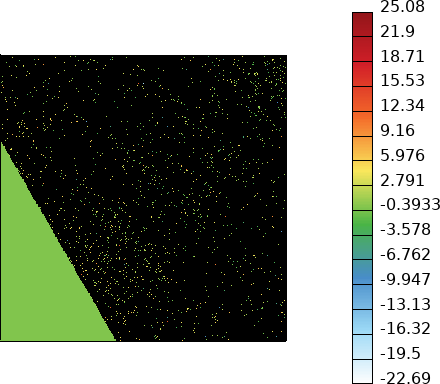}
	\end{minipage}
	\begin{minipage}[t]{0.33\textwidth}
		\centering (c)  $\beta_h=\max\left\{1,\frac{|\mathbf v| h}{2\mu}\right\}$\vskip0.5cm
		
		\includegraphics[width=0.9\textwidth]{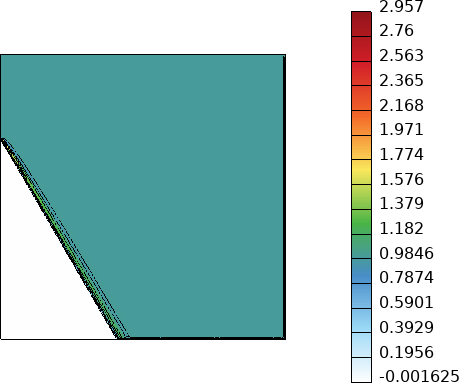}
	\end{minipage}
	\vskip0.5cm
	
	\begin{minipage}[t]{0.33\textwidth}
		\includegraphics[width=0.9\textwidth]{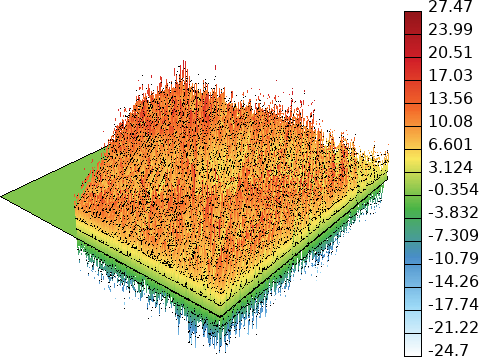}
	\end{minipage}%
	\begin{minipage}[t]{0.33\textwidth}
		\includegraphics[width=0.9\textwidth]{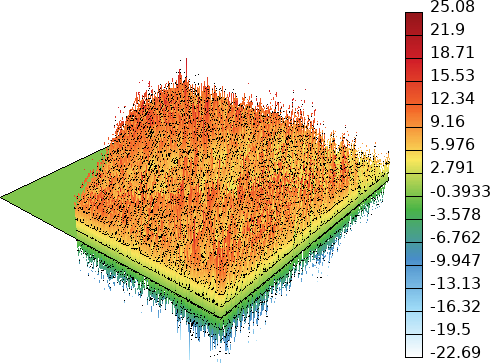}
	\end{minipage}
	\begin{minipage}[t]{0.33\textwidth}
		\includegraphics[width=0.9\textwidth]{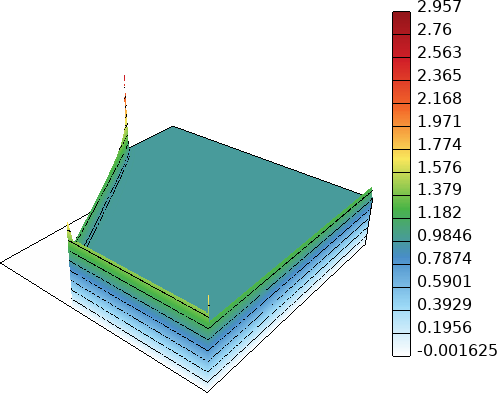}
	\end{minipage}
	\vskip0.25cm 
	\caption{Numerical solutions to the stationary interface convection-diffusion problem (uniform mesh, $h=\frac{1}{512}$, $\delta = 0$). First row: top views. Second row: elevation is proportional to the value of the plotted function.}
	\label{fig:kd}
\end{figure}

In Fig.~\ref{fig:kd},
we show the numerical solutions calculated on the uniform mesh with
spacing $h=\frac{1}{1024}$. It can be seen that failure to properly
stabilize the discretized
convective terms results in global spurious oscillations. The
\rev{projected gradient} stabilization term with the proposed choice of the parameter
$\beta_h$ mitigates this effect and localizes violations of the
discrete maximum principle \rev{(DMP)} to layers of elements with
steep gradients without the need for additional modifications of
the discretized weak form. \rev{Moreover, the magnitude of spurious
undershoots
  and overshoots is significantly reduced (see the 3D plots in
  Fig.~\ref{fig:kd}). If violations of the DMP are unacceptable, they can be avoided, e.g., using the framework of algebraic flux correction \cite{AFC,FCTbook,WSbook}.}

\section{Conclusions}
\label{sec:conclusions}
The paper studied an unfitted FEM with \rev{\rev{projected gradient}} stabilization. Two versions of the method, the sharp interface and diffuse interface approaches, were introduced and applied to numerically solve an elliptic interface problem. The \rev{\rev{projected gradient}} stabilization was applied globally, rather than being restricted to a narrow strip around the interface.

A priori analysis of the \rev{\rev{projected gradient}} stabilization was performed for the case of $P^1$
finite elements. The resulting sharp interface version of the FEM was proven to be optimal-order convergent. The method's performance was compared to the unstabilized variant from \cite{hansbo2}. The proposed version demonstrated comparable accuracy while being algebraically stable, providing an extension of the solution beyond the physical domain and enabling a diffuse interface implementation.
The diffuse interface variant showed similar accuracy but currently lacks numerical analysis.

We conclude that \rev{projected gradient} stabilization is a simple and effective technique that facilitates the straightforward implementation of unfitted FEM. Handling higher-order elements remains an open problem for future research.

\section*{Acknowledgments} M.O. was supported in part by the U.S. National Science Foundation under awards DMS-2309197 and DMS2408978. 

\end{document}